\documentclass[12pt,a4paper]{amsart}

\usepackage[american]{babel}
\usepackage{bm,amsmath,amssymb}
\usepackage{amsbsy,amsfonts,array}
\usepackage{verbatim,graphicx}
\usepackage{epstopdf}
\usepackage{subfigure}
\graphicspath{{./figs/},{./}}
\usepackage[margin=3cm]{geometry} 

\newcommand{\Div}{\text{div}}
\newcommand{\nabu}{\boldsymbol{\nabla}}
\newcommand{\Frac}[2]{\displaystyle\frac{#1}{#2}}
\newcommand{\Sum}[2]{\displaystyle{\sum\limits_{#1}^{#2}}}
\newcommand{\vect}[1]{\bm{#1}}

\newtheorem{theorem}{Theorem}

\title[Photocurrent transients in OSCs]{Analytical and Numerical Study of Photocurrent Transients
in Organic Polymer Solar Cells}
\address{Dipartimento di Matematica \lq\lq F.Brioschi\rq\rq,
               Politecnico di Milano, \\
	             via Bonardi 9, 20133 Milano Italy}
\author[C. de Falco]{Carlo de Falco}
\author[R. Sacco]{Riccardo Sacco}
\author[M. Verri]{Maurizio Verri}

\begin{document}
\maketitle

\begin{abstract}
This article
is an attempt to provide a self consistent picture, including 
existence analysis and numerical solution algorithms, 
of the mathematical problems 
arising from modeling photocurrent transients in Organic-polymer 
Solar Cells (OSCs).
The mathematical model for OSCs consists of a
system of nonlinear diffusion-reaction partial 
differential equations (PDEs) with electrostatic convection,
coupled to a kinetic ordinary differential equation (ODE).
We propose a suitable reformulation of the model
that allows us 
to prove the existence 
of a solution 
in both stationary and transient conditions 
and
to better highlight 
the role of exciton dynamics in determining the device turn-on time. 
For the numerical treatment of the problem, we 
carry out a temporal semi-discretization using an implicit 
adaptive method, and the resulting sequence of 
differential subproblems is linearized using the Newton-Raphson
method with inexact evaluation of the Jacobian. 
Then, we use exponentially fitted finite elements for the
spatial discretization, and we carry out a thorough validation
of the computational model by extensively investigating 
the impact of the model parameters on photocurrent transient times.
\end{abstract}

\section{Introduction and Motivation}\label{sec:intro}
A continuously growing attention has been paid
over the last years by the international community and government 
authorities to monitoring the effect of the increase of global
concentrations of carbon dioxide, methane and nitrous oxide
on the quality of our everyday life. The results of the investigation 
carried out by the Intergovernmental Panel on Climate Change~\cite{IPCC2007}
have brought the European Union (EU) to the decision that
carbon dioxide emissions should decrease by 20 percent,
and that 20 percent of the energy produced in EU should originate 
from renewable energy sources, such as wind, water, biomass, and solar, 
not later than 2020~\cite{EUCouncilMeeting2007}.
In this perspective, research and design 
of third generation (3G) photovoltaic devices~\cite{GreenBook} 
for solar energy conversion into electrical and thermal energy 
turns out to be a central topic in the wider area of renewable energy sources. 
Roughly speaking, 3G photovoltaic devices 
can be divided into two main classes: 
electrochemical cells~\cite{Gratzel:2001eu,Graetzel2005,graetzel2008}
and organic-polymer cells~\cite{Mihailetchi2004,McGehee2007,Gunes2007} which are the topic of the present article. 
Most of investigation activity in solar cell design is devoted to the experimental study 
of innovative materials 
for efficient and flexible 
technologies, and is not presently accompanied by a 
systematic use of computational models 
to predict and optimize their performance. 
This article is an attempt to fill this gap by introducing the numerical
engineering community to the mathematical problems that arise
in the context of modeling and simulation of OSCs. With this aim, we try  to provide
a reasonably self-contained picture of the topic, including a discussion of the peculiarities of the model, an analysis of the existence of a solution, and the description of a robust computational algorithm to compute such solution.
In particular, we focus on 
a special class of OSCs, namely that of Bulk Hetero-Junction (BHJ) devices, 
that currently  represent the most promising 
technology in terms of energy conversion efficiency~\cite{Gunes2007,McGehee2007}.
Charge transport in BHJs is described 
by a set of nonlinear PDEs of diffusion-reaction type
with electrostatic convection coupled with a kinetic ODE
for the temporal evolution of exciton concentration
in the cell~\cite{Barker2003,Buxton2007,Hwang2008,coehoorn}. 
Sect.~\ref{sec:physics} is devoted to the description of the structure and working principles of BHJs while  in Sect.~\ref{sec:model} the mathematical model is introduced and the connection between its features and the physical phenomena involved in photocurrent generation is drawn. Some effort is also put into highlighting the main differences between the problem at hand and
the case of more standard crystalline inorganic semiconductor devices.
In Sect.~\ref{sec:analysis}, under suitable assumptions on the
model coefficients, {\it i}) we prove the existence
of a solution of the problem in stationary conditions;
and {\it ii}) we derive a simplified model in transient
conditions, that is amenable for a qualitative analysis
of the time response of the device, and for which we again 
prove existence of a solution.
For the numerical treatment of the problem, which is the topic of Sect.~\ref{sec:numericalmethod}, we 
carry out a temporal semi-discretization using an implicit 
adaptive method, and the resulting sequence of 
differential subproblems is linearized using the Newton-Raphson
method with inexact evaluation of the Jacobian. 
Then, we use exponentially fitted finite elements for the
spatial discretization, to ensure a stable
approximation of the internal and boundary layers
arising in the distribution profile of the
photogenerated carriers.
The numerical experiments of Sect.~\ref{sec:numericalresults} are meant, on the one hand, to illustrate the complex interplay among different physical phenomena determining the photocurrent turn-on transient time of a realistic BHJ cell in different regimes and, on the other hand, to characterize the range of applicability of the reduced model introduced in Sect.~\ref{sec:analysis}. 
In Sect.~\ref{sec:conclusions} we address some concluding remarks 
and indicate possible future research directions.

\section{Bulk Heterojunction Organic Solar Cells}\label{sec:physics}
\begin{figure}
\begin{center}
\includegraphics[width=.6\linewidth]{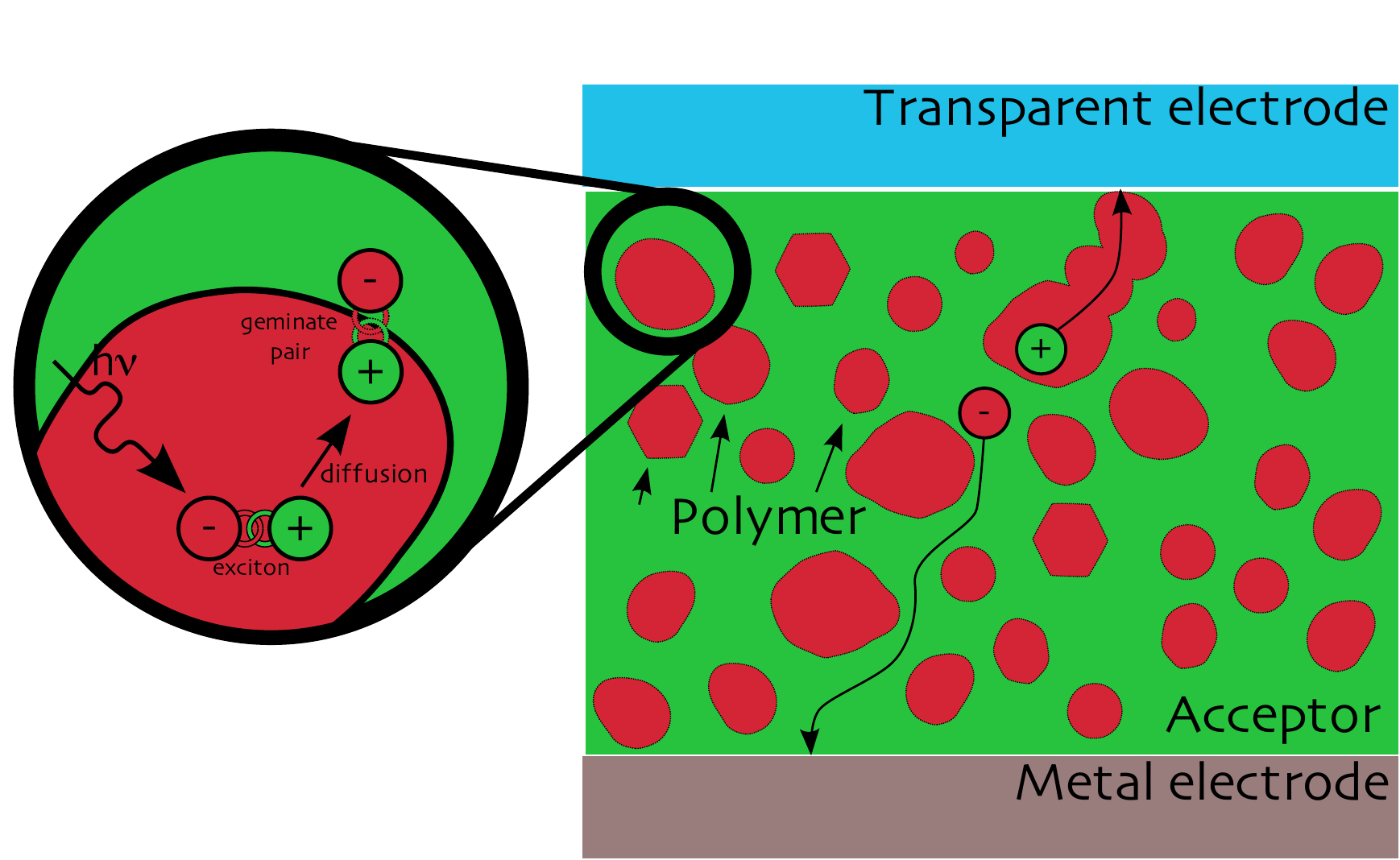}
\caption{Bulk Heterojunction OSCs.}
\label{fig:bhj}
\end{center}
\end{figure}

Before presenting the mathematical model which is the main focus of this paper, a schematic 
description of working principle of OSCs, and in particular of those with
a BHJ structure, is in order.
For more details on the subject the interested reader is referred to~\cite{Gunes2007,McGehee2007}.
The simplest possible structure for an organic-polymer based solar cell is depicted in Fig.~\ref{fig:planarosc}:
\begin{figure}
\begin{center}
\subfigure[]{\includegraphics[width=.32\linewidth]{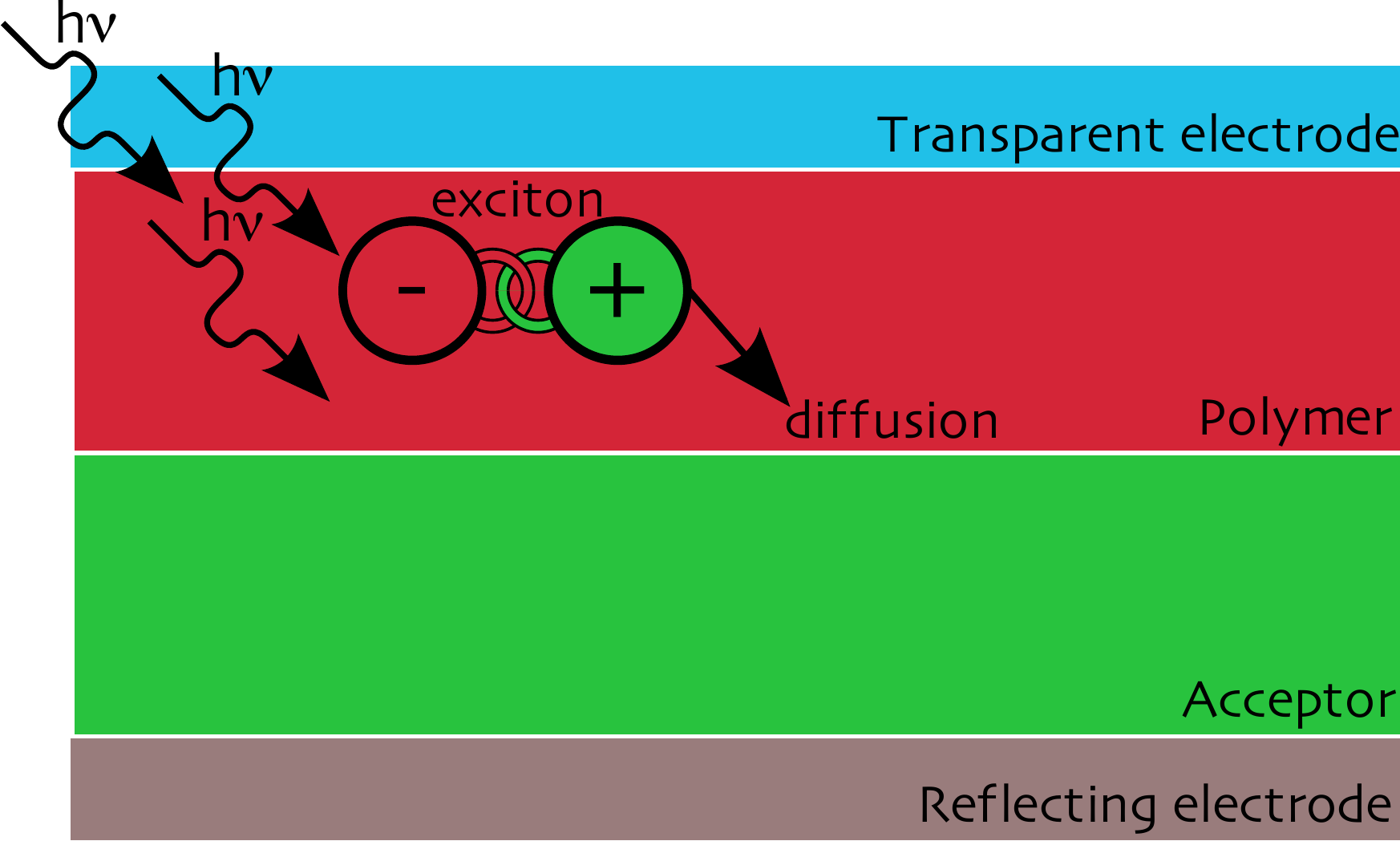}\label{fig:exciton}}
\subfigure[]{\includegraphics[width=.32\linewidth]{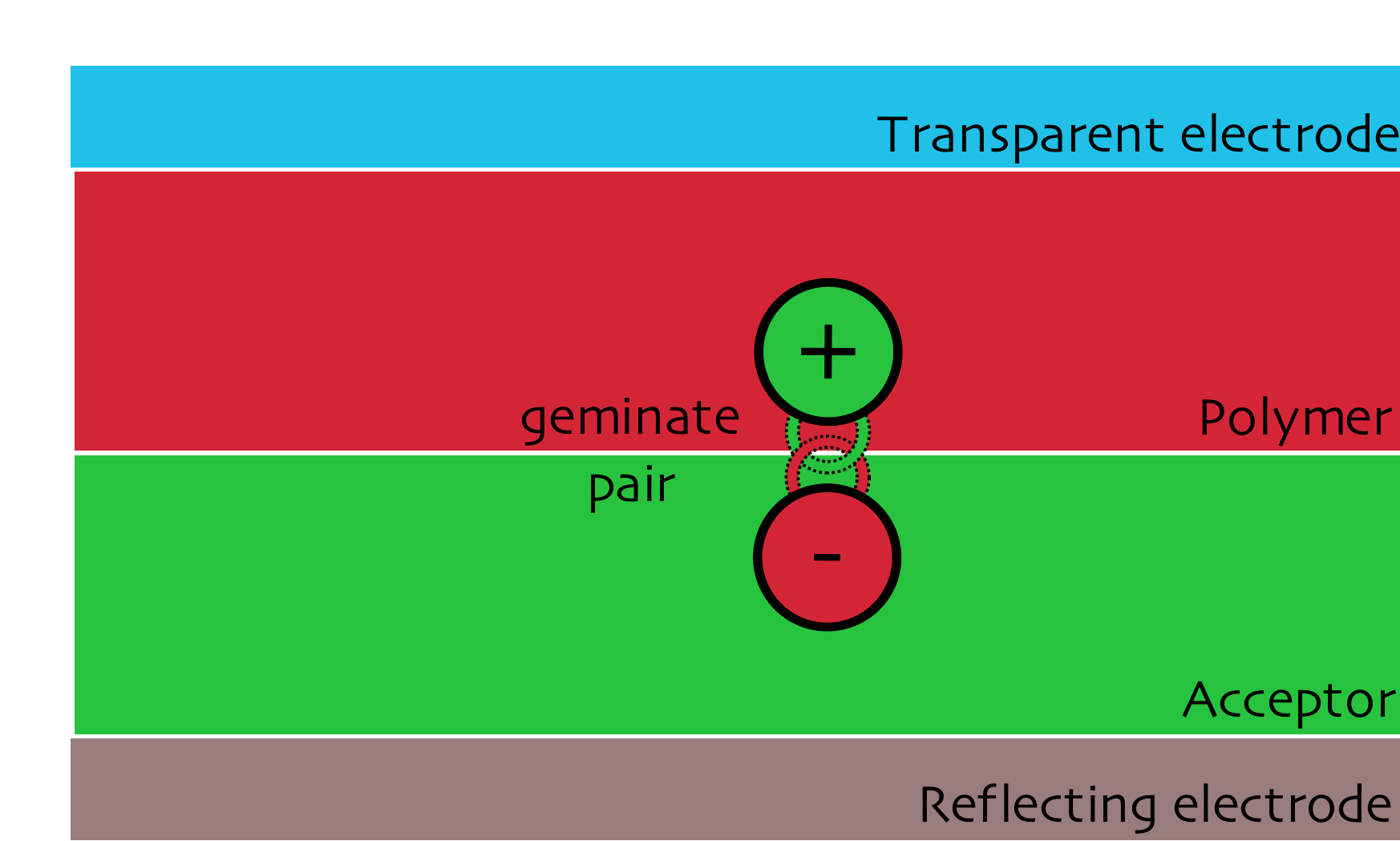}\label{fig:geminatepair}}
\subfigure[]{\includegraphics[width=.32\linewidth]{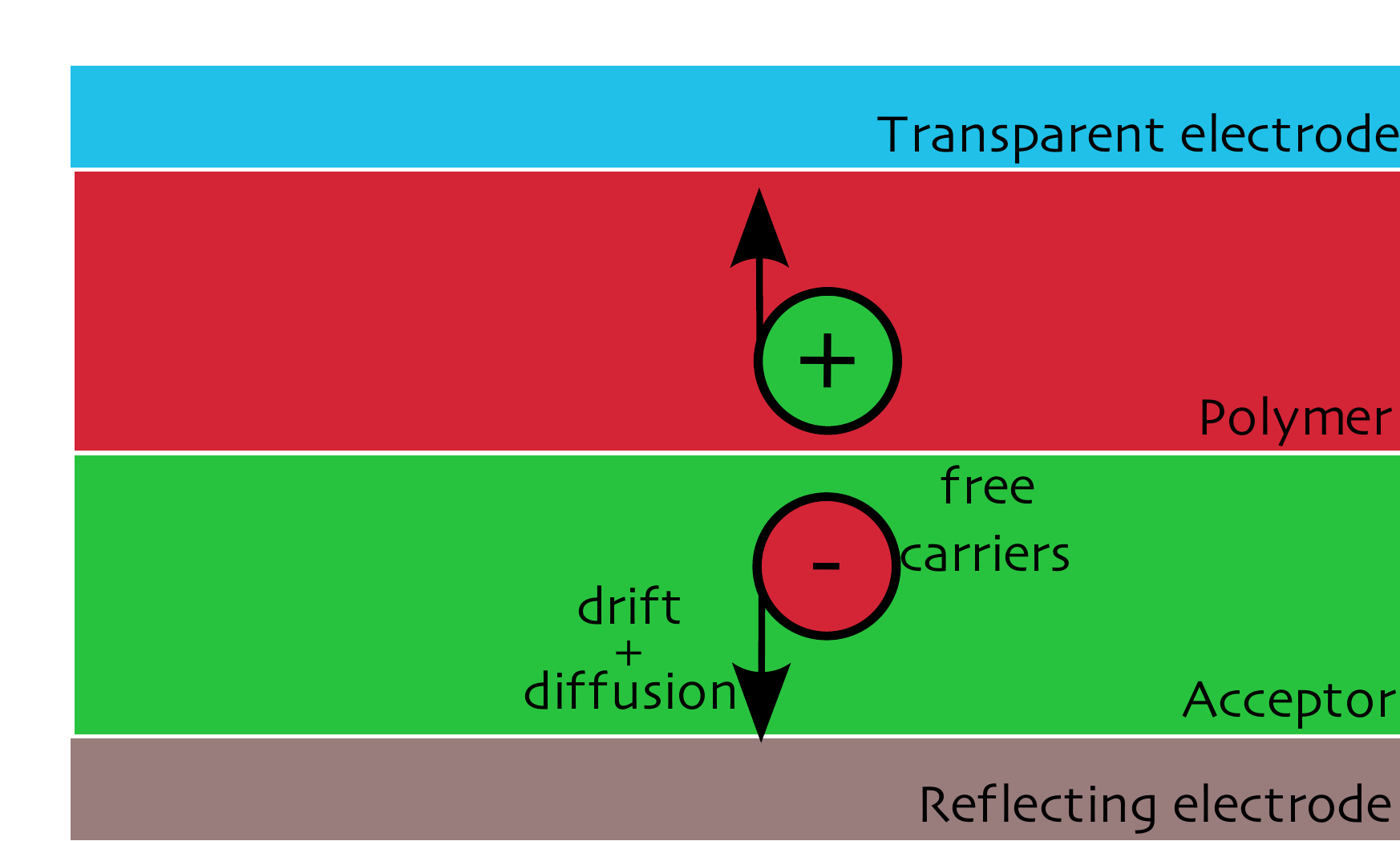}\label{fig:harvesting}}
\caption{Working principle of OSCs.}
\label{fig:planarosc}
\end{center}
\end{figure}
two thin films composed of a conjugated organic polymer and of a material with high electron affinity, usually referred to as a \emph{acceptor}
 are sandwiched between one transparent ({\it e.g.} indium-tin-oxide or fluorinated tin oxide) and one reflecting metal contact (usually aluminum or silver).
When illuminated, electrons in the Highest Occupied Molecular Orbital (HOMO) in the polymer are promoted to
the Lowest Unoccupied Molecular Orbital (LUMO) thus forming an electron-hole pair.
Such pair, which we refer to as an \emph{exciton} (Fig.~\ref{fig:exciton}), in contrast to what is usually the case in standard inorganic semiconductors, 
is electrically neutral and has very strong binding energy (of the order of $1$eV) with a radius in the sub-nanometer range.
The \emph{diffusion length} $\lambda_{X}$ 
of a moving exciton in commonly used polymeric materials is of the order of a few nanometers. An exciton has a non-negligible chance of eventually reaching the polymer/acceptor interface only if it was  photo-generated within a distance $\leq \lambda_{X}$. In case this occurs, the
built-in chemical potential drop produced by the difference in electron affinity between the two materials
is strong enough to \emph{stretch} the exciton driving  the electron and hole to a distance of the order of $1$nm thus reducing the strength of their Coulomb attraction. This less tightly bound electron-hole pair is referred to in the literature as a \emph{geminate pair} (Fig.~\ref{fig:geminatepair}) and the energy of the bond is low enough that it can be overcome by the electric field induced by a small voltage difference applied at the contacts. 
The newly separated electron and hole migrate, driven by electric field drift and diffusion forces, to the anode and cathode, respectively, where they are \emph{harvested} thus producing a net current (Fig.~\ref{fig:harvesting}).
The currently investigated most promising device technology to maximize the efficiency of the photogeneration process is 
the BHJ cell depicted in Fig.~\ref{fig:bhj} which is produced by spin-casting both the polymer (usually rr-P3HT or MDMO-PPV) and the acceptor (usually some derivative of fullerene or inorganic nanoparticles, {\it e.g.} titanium-dioxide) from a common solution. This process results in a highly folded structure that has the advantage that all photo-generated excitons  eventually reach an interface, at the price of reducing the \emph{effective} carrier mobility 
because of the convoluted path that carriers need to travel to reach the contacts.
Also, from a perspective that is more relevant to the topic of this paper, the highly disordered structure of BHJs makes it difficult to characterize model parameters, as an averaging over the highly disordered nanostructure of the device would be required. Therefore the typical approach is to estimate the parameter values experimentally and resort to numerical simulations to properly interpret the measurement results.

\section{The Mathematical Model}\label{sec:model}
In this section we illustrate the mathematical
model of the photogeneration mechanisms that drive charge 
transport in BHJ solar cells (see \cite{Barker2003,Gunes2007,Buxton2007,Hwang2008,coehoorn}).
The polymer/acceptor blend is represented by a homogeneous material 
filling  a bounded domain $\Omega \subset \mathbb{R}^d$, $d\geq 1$,
with a Lipschitz boundary $\Gamma \equiv \partial \Omega$ divided 
into two disjoint subregions, $\Gamma_D$ and $\Gamma_N$, representing
the interface between metal and polymer blend and 
interior artificial boundaries, respectively.
We assume that ${\rm meas\ }(\Gamma_D)>0$ and $\Gamma_D \cap 
\Gamma_N = \emptyset$, and denote by 
$\boldsymbol{\nu}$ the outward unit normal vector along $\Gamma$.
\subsection{Governing Equations}
Charge transport in the device is governed by the set of continuity 
equations
\begin{subequations}\label{eq:modelequations}
\begin{equation}\label{eq:continuity}
\left\{
\begin{array}{llll}
\Frac{\partial n}{\partial t} & - \Div \vect{J}_{n} & = & 
G_n - R_n\ n\\[2mm]
\Frac{\partial p}{\partial t} & - \Div \vect{J}_{p} & = & 
G_p - R_p\ p
\end{array}
\right.\qquad \mbox{in }\Omega_{T},
\end{equation}
where $\Omega_{T}\equiv \Omega \times (0,T)$, $T > 0$, $n$ and $p$ denote the \emph{electron and hole density}, respectively. Using from now on the symbol $\eta$ to indicate either of $n$ or $p$,
$\vect{J}_{\eta}$ are the corresponding \emph{flux densities}, $G_{\eta}$, are the carrier \emph{generation rates}, and $R_{\eta} \eta$ are the \emph{recombination rates}. 
As electrons are negatively charged while hole charge is positive, the 
\emph{total current density} $\vect{J}$ can be expressed as 
$\vect{J} = q\left( \vect{J}_{p} - \vect{J}_{n}\right)$ 
where $q>0$ is the magnitude of the electron charge.
The charge carrier flux densities are, in turn, 
each composed of an \emph{electrostatic drift} term 
and a \emph{diffusion} term
\begin{equation}\label{eq:driftdiffusion}
\left\{
\begin{array}{ll}
\vect{J}_n & = D_n
\nabu n - \mu_{n} \, n \, \nabu \varphi\\[2mm]
\vect{J}_p & = D_p
\nabu p + \mu_{p} \, p \, \nabu \varphi \\[2mm]
\end{array}
\right.\qquad \mbox{in }\Omega_{T},
\end{equation}
$D_{\eta}$ being the charge carrier \emph{diffusion coefficients} and $\mu_{\eta}$ the \emph{carrier mobilities}. 
The \emph{electrostatic potential} $\varphi$ satisfies the Poisson 
equation
\begin{equation}\label{eq:poisson}
- \Div(\varepsilon \nabu \varphi) = q (p - n) \qquad \mbox{in }\Omega_{T}, 
\end{equation}
where $\varepsilon$ is the (averaged) dielectric permittivity of the blend. Notice that, as there are 
usually no dopants in organic cells, the net charge density on the right-hand-side in~\eqref{eq:poisson} is given by the carrier densities only.
We denote by $X$ the \emph{volume density of geminate pairs} and we express its rate of change as
\begin{equation}\label{eq:kinetic}
\displaystyle \frac{\partial X}{\partial t} = g - r\qquad \mbox{in }\Omega_{T}.
\end{equation}
\end{subequations}
The geminate-pair generation rate $g$ in~\eqref{eq:kinetic} 
can be split into two contributions as
\begin{equation}\label{eq:excitongeneration}
g = \underbrace{G(\vect{x}, t)}_{\mbox {(\it a)}} + 
\underbrace {\gamma \, p \, n}_{\mbox{(\it b)}},
\end{equation}
({\it a}) accounting for the  rate at which excitons reach the material interfaces and are
partially separated and ({\it b}) accounting for 
the rate at which free electrons and holes are attracted to each other and recombine. Process ({\it b}) is referred to as
\emph{bimolecular recombination} and the coefficient $\gamma$ is described according to the Langevin theory~\cite{Mihailetchi2004}. The rate of process ({\it a}) is equal to the rate $G(\vect{x}, t)$ at which photons are absorbed, which we assume in what follows to be a known function of  position and time.
As in BHJ all excitons are eventually transformed in geminate pairs it is legitimate, with a slight abuse of notation, to use in the following the two terms as synonyms. 
As for the term $r$ in~\eqref{eq:kinetic} it can also be split into two contributions as
\begin{equation}\label{eq:excitonrecombination}
r = \underbrace{k_{diss} X}_{\mbox {(\it c)}} + \underbrace {k_{rec} X}_{\mbox{(\it d)}},
\end{equation}
{\it (c)} accounting for the rate at which geminate pairs that are not split recombine and {\it (d)} accounting for the rate at which free electrons and holes are produced by separation of a bound pair. We assume the coefficient $k_{rec}$ to be a given constant while $k_{diss}$ depends on the magnitude of the electric field $\vect{E} = -\nabu \varphi$ 
as described in~\cite{Mihailetchi2004}.
As we assume free carriers to be generated only by dissociation of a geminate pair and to be annihilated only by recombination into a geminate pair, the generation rates satisfy $G_{n} = G_{p} = k_{diss} X$ while for the recombination rates $R_{n} n = R_{p} p = \gamma p n$ holds.

We wish at this point to stress some peculiarities of the model we have introduced compared to the standard case of crystalline inorganic semiconductors. 
The main difference is represented by the strong influence that the exciton reaction kinetics described by equation~\eqref{eq:kinetic} has on device performance. Indeed, such a kinetics affects both the energy conversion efficiency in the steady state operation and the  turn-on transient time.  This latter, in particular, is relevant for the characterization of material properties that can not be determined by first-principles because of the highly convoluted device nanostructure.
Furthermore, although equations~\eqref{eq:continuity}-\eqref{eq:driftdiffusion} are analogous to those describing 
charge transport in ordered inorganic semiconductors, the physical driving mechanisms at the microscopic level  are quite different. In particular, while in monocrystalline semiconductors charge carriers are essentially free to move within \emph{delocalized} orbitals, in the materials we study here transport happens via \emph{hopping} of charges between \emph{localized} orbitals. 
This microscopic difference is reflected in the macroscopic models for the diffusion and mobility coefficients for organic semiconductor materials which {\it (i)} introduce very different dependencies on temperature and electric field magnitude~\cite{gill1972,horowitz1998}, and {\it (ii)}  introduce a dependency of the mobility on the carrier densities~\cite{coehoorn}. 
\subsection{Boundary and Initial Conditions}
A delicate and important issue is that of devising 
a set of boundary conditions to accurately describe the complex
phenomena of charge injection and recombination occurring at the
interface $\Gamma_D$ separating the metal contacts from the 
semiconductor bulk. Precisely, according to~\cite{CampbellScott1999,Barker2003},
such conditions 
are expressed in the following Robin-type form
\begin{subequations}\label{eq:boundaryconditions}
\begin{align}\label{eq:electronboundaryconditions}
\kappa_n \, \vect{J}_n \cdot \boldsymbol{\nu} =
\beta_n - \alpha_n \, n \qquad \mbox{on } \Gamma_D\times (0,T)\\[2mm]
\label{eq:holeboundaryconditions}\kappa_p \, \vect{J}_p \cdot \boldsymbol{\nu} =
\beta_p - \alpha_p \, p \qquad \mbox{on } \Gamma_D\times (0,T),
\end{align}
where $\kappa_\eta$ are non negative parameters while
$\beta_\eta$ are the rates at which charges are injected
into the device and $\alpha_\eta \eta$ are the rates at which 
electrons and holes recombine with their image charges at the contacts, respectively.
Reliable models for the above parameters are
still subject of extensive investigation as the basic description
proposed in the milestone reference~\cite{CampbellScott1999}
needs to be modified via empirical fitting to avoid the 
occurrence of unphysical behavior in the computed 
solution~\cite{Hwang2008,Hwang2009_private}.
As for the electric potential, the Dirichlet condition
\begin{equation}\label{eq:potentialboundaryconditions}
\varphi = \Psi_{D}\qquad \mbox{on } \Gamma_D\times (0,T)
\end{equation}
is enforced, where the datum $\Psi_{D}$ accounts for both the externally applied voltage and the work-function difference between the contact materials.
On $\Gamma_N$, which represents the interior artificial boundary, homogeneous Neumann conditions for the flux densities and the electric field are imposed.
Finally, positive initial conditions $n(\vect{x}, 0)=n_{0}(\vect{x})$, $p(\vect{x}, 0)=p_{0}(\vect{x})$, and $X(\vect{x}, 0)=X_{0}(\vect{x})$ are needed to complete the mathematical model.
\end{subequations}

\section{System Analysis of the Model}\label{sec:analysis}
In this section, we deal with the analysis of the existence
of a solution of system~\eqref{eq:continuity}-\eqref{eq:kinetic}
in both stationary and transient regimes, under the following assumptions:
\begin{description}
\item[(H1)] $\gamma$, 
$k_{diss}$, $k_{rec}$ and $G$ are all positive constant 
quantities in $\Omega_T$;
\item[(H2)] $D_{\eta} = V_{th} \mu_{\eta}$,
$V_{th}$ being the \emph{thermal voltage} and 
$\mu_{\eta} \ge \mu_{\eta_{0}} > 0$ a.e. in $\Omega_{T}$;
\item[(H3)] $v_{n}, v_{p} \le v^{max} < + \infty$ where $v_{\eta} := \mu_{\eta} |\vect{E}|$;
\item[(H4)] $\kappa_{\eta}=0$ and $\alpha_{\eta}, \beta_{\eta}$ are functions of position only in~\eqref{eq:electronboundaryconditions}-\eqref{eq:holeboundaryconditions}.
\end{description}

Although the purpose of the set of hypotheses (H1)-(H4) is mainly to reduce the 
mathematical complexity of the problem, we wish here to comment about their physical 
plausibility.
Assumption (H1) allows us to express in
an easy manner the dependent variable $X$ as a function of
$n$, $p$ and of the input data $G$ and $X_0$, in such a way that the 
resulting equivalent system (in the reduced set of 
unknowns $\varphi$, $n$ and $p$) can be
written in the form of a two-carrier drift-diffusion (DD) model.
As the coefficients involved in (H1) depend, in general, only on the magnitude of the 
electric field,
such an assumption is reasonable if the field itself varies weakly within 
the simulation domain, 
which is often the case in realistic photovoltaic devices
as is confirmed by the numerical experiments of Sect.~\ref{sec:numericalresults}.
Assumption (H2) is the classical Einstein relation valid in inorganic semiconductors 
and corresponds to neglecting the (higher order) effect of energetic disorder~\cite{coehoorn}.
The saturation of convective velocities expressed by assumption (H3) is reasonable in a 
structure that is highly folded as that of BHJs and is indeed commonly employed in commercial packages for organic semiconductor simulation~\cite{atlas}.
Assumption (H4) corresponds to an infinite carrier recombination rate at the contacts.

\subsection{Stationary Regime}\label{sec:stationary-analysis}
Setting $\partial X /\partial t = 0$ in~\eqref{eq:kinetic},
we can eliminate the dependent variable $X$ in favor of 
$n$, $p$ and of the input function $G$, to obtain
\begin{equation}\label{eq:x} 
X(\vect{x})  =  \tau G +  \gamma \tau  p(\vect{x}) n(\vect{x}) 
\end{equation}
where 
\begin{equation}\label{eq:tau}
\tau := \displaystyle\frac{1}{k_{diss}+k_{rec}}
\end{equation}
is the time of response of the 
generation/recombination terms to light stimuli.
Using~\eqref{eq:x}-\eqref{eq:tau} and (H4), the stationary OSC model reads:
\begin{equation}\label{eq:model-stationary}
\left\{
\begin{array}{lll}
-\Div(\varepsilon \nabu \varphi) & = & q (p - n) \\[1mm]
-\Div \vect{J}_n  & = & 
\tau \left(k_{diss} \; G -  \gamma \; k_{rec} \; p n\right) \\[1mm]
-\Div \vect{J}_p & = & 
\tau \left(k_{diss} \; G -  \gamma \; k_{rec} \; p n\right),
\end{array}
\right.\qquad \mbox{in }\Omega
\end{equation}
supplied with the boundary conditions
\begin{equation}\label{eq:bcs-2}
\left\{
\begin{array}{ll}
\varphi=\Psi_D, \quad n = n_{D} := \Frac{\beta_{n}}{\alpha_{n}}, \quad
p = p_D  := \Frac{\beta_{p}}{\alpha_{p}} & \quad \mbox{on } \Gamma_D \\[1mm]
\vect{J}_n \cdot \boldsymbol{\nu} = \vect{J}_p \cdot \boldsymbol{\nu}
= \nabu \varphi \cdot \boldsymbol{\nu} = 0 & 
\quad \mbox{on } \Gamma_N.
\end{array}
\right.
\end{equation}
\begin{theorem}[Existence of a solution in stationary regime]
\label{th:steady_state_ex}
Let assumptions (H1)-(H4) be satisfied
and $(\Psi_D, n_D, p_D) \in (L^{\infty}(\Gamma_D))^3$.
Then, problem~\eqref{eq:model-stationary}--\eqref{eq:bcs-2}
admits a weak solution $(\varphi^{\ast}, u^{\ast}, v^{\ast})
\in (H^1(\Omega) \cap L^{\infty}(\Omega))^3$ and there exist 
positive constants $\underline{\mathcal{M}}$, $\overline{\mathcal{M}}$,
$\underline{\mathcal{K}}$, $\overline{\mathcal{K}}$
such that
\begin{equation}\label{eq:estimate-stationary}
\underline{\mathcal{M}}\ \leq n^{\ast}, \; 
p^{\ast} \leq \overline{\mathcal{M}},
\qquad
\underline{\mathcal{K}}\ \leq \varphi^{\ast} \leq \overline{\mathcal{K}}
\qquad
 \mbox{{\rm a.e. }in} \;  \Omega.
\end{equation}
\end{theorem}
The proof of Theorem~\ref{th:steady_state_ex}
follows closely the guidelines of~\cite{Mar90}, Sect.3.3 and is sketched below.
Using (H2) we can write the two flux densities as  
\begin{equation}\label{eq:fluxes-slotboom}
\left\{
\begin{array}{lll}
\vect{J}_n & = & \mu_n V_{th} \; n_r \; \mathrm{e}^{\varphi/V_{th}} \nabu u, \\[2mm]
\vect{J}_p & = & \mu_p V_{th} \; n_r \; \mathrm{e}^{-\varphi/V_{th}} \nabu v, 
\end{array}
\right.
\end{equation}
where the new (dimensionless) dependent variables $u$ and $v$ are related to
the carrier densities $n$ and $p$ by the Maxwell--Boltzmann statistics
\begin{equation}\label{eq:Maxw-Boltz}
n = n_r \; u \mathrm{e}^{\varphi/V_{th}}, \qquad 
p = n_r \; v \mathrm{e}^{-\varphi/V_{th}},
\end{equation}
$n_r >0$ being a reference 
concentration. 
System~\eqref{eq:model-stationary}-\eqref{eq:bcs-2} then becomes:
\begin{equation}\label{eq:model-stationary-2}
\left\{
\begin{array}{lll}
-\Div(\varepsilon \nabu \varphi) & = & 
q \; n_r (u \mathrm{e}^{\varphi/V_{th}} - v \mathrm{e}^{-\varphi/V_{th}})\\[2mm]
-\Div (\mu_n V_{th}  \; \mathrm{e}^{\varphi/V_{th}} \nabu u) & = & 
\Frac{\tau k_{diss} G}{n_r} (1 - uv) \\[3mm] 
-\Div (\mu_p V_{th} \; \mathrm{e}^{-\varphi/V_{th}} \nabu v) & = & 
\Frac{\tau k_{diss} G}{n_r} (1 - uv)
\end{array}
\right.\qquad \mbox{in }\Omega
\end{equation}
and
\begin{equation}\label{eq:bcs-3}
\left\{
\begin{array}{ll}
\varphi=\Psi_D, \;\; u = u_{D}:=\Frac{n_D}{n_r} \; \mathrm{e}^{-\Psi_D/V_{th}}, \;\;
v =v_{D}:= \Frac{p_D}{n_r} \; \mathrm{e}^{\Psi_D/V_{th}} & \quad \mbox{on } \Gamma_D \\[2mm]
\vect{J}_n \cdot \boldsymbol{\nu} = \vect{J}_p \cdot \boldsymbol{\nu}
= \nabu \varphi \cdot \boldsymbol{\nu} = 0 & 
\quad \mbox{on } \Gamma_N. 
\end{array}
\right.
\end{equation}
Using the boundedness of the Dirichlet data, the positivity of $n_{D}$ and $p_{D}$ and
choosing $n_r$ in such a way that $(\gamma k_{rec} n_r^2)/(k_{diss} G)=1$,
we can see that 
\begin{equation}\label{eq:a-priori-bound}
\mathrm{e}^{-\Psi^+/V_{th}} \leq u_D, \; v_D
\leq \mathrm{e}^{\Psi^+/V_{th}} 
\qquad \qquad \mbox{a.e. on} \;  \Gamma_D,
\end{equation}
where
$$
\Psi^+:= \max \left\{ \max (\displaystyle 
\sup_{\Gamma_D} (-\varphi_{nD}), \sup_{\Gamma_D} (\varphi_{pD})), 
-\min(\displaystyle 
\inf_{\Gamma_D} (-\varphi_{nD}), \inf_{\Gamma_D} (\varphi_{pD})) \right\}
$$
and
$$
\varphi_{nD}:= \Psi_D - V_{th} \ln(n_D/n_r),
\qquad  
\varphi_{pD}:= \Psi_D + V_{th} \ln(p_D/n_r).
$$
Then, by applying Theorem 3.3.16 of~\cite{Mar90} to system~\eqref{eq:model-stationary-2}-\eqref{eq:bcs-3} and going back to the original variables $n$ and $p$ via the inversion of~\eqref{eq:Maxw-Boltz}, we conclude that Theorem~\ref{th:steady_state_ex} holds with
\begin{subequations}
\begin{align}
\underline{\mathcal{K}} &= n_r \; \mathrm{e}^{\displaystyle -\widehat{\Psi}^+/V_{th}}, &\quad
\overline{\mathcal{K}}  &= n_r \; \mathrm{e}^{\displaystyle \widehat{\Psi}^+/V_{th}} \\
\underline{\mathcal{M}} &= \min \left(\displaystyle \inf_{\Gamma_D} \Psi_D, \; - \Psi^+ \right), &\quad
\overline{\mathcal{M}}  &= \max \left( \displaystyle \sup_{\Gamma_D} \Psi_D, \; \Psi^+ 
\right)
\end{align}
\end{subequations}
where $\widehat{\Psi}^+:=\displaystyle \sup_{\Gamma_D} |\Psi_D| + \Psi^+$.

\subsection{Transient Regime}\label{sec:transient-analysis}
Analogously to what we have done in Sect.~\ref{sec:stationary-analysis} 
in the stationary case, we can use~\eqref{eq:kinetic}
to eliminate the dependent variable $X$ in favor of $n$, $p$ and 
of the input functions $G$ and $X_0$, to obtain
\begin{eqnarray}\label{eq:x_di_t} 
X(\vect{x},t) & = & 
\xi(\vect{x},t) + \gamma \displaystyle 
\int_{0}^{t } p(\vect{x}, s) \; n(\vect{x}, s) 
\mathrm{e}^{-(t-s)/\tau} \ ds,
\end{eqnarray}
where $\xi(\vect{x},t) := X_{0}(\vect{x}) 
\mathrm{e}^{\displaystyle -t/\tau} + \tau \; G (1 - \mathrm{e}^{\displaystyle -t/\tau})$.
The quadratic convolution term in~\eqref{eq:x_di_t} makes the dependence 
of the current on the electron and hole densities non-local 
in time with a ``memory window'' of size proportional to $\tau$.
For the subsequent existence analysis it is convenient
 to manipulate such term so that
we can write the 
continuity equations in the following equivalent form:
\begin{equation}\label{eq:ncont_no_x_2}
\left\{
\begin{array}{l}
\Frac{\partial n}{\partial t} - \Div \vect{J}_{n} =
k_{diss} \xi- \gamma \, \tau 
\left( k_{rec} + k_{diss} \mathrm{e}^{-t/\tau}\right) 
\, p \; n + \ I \\[.3cm]
\Frac{\partial p}{\partial t} - \Div \vect{J}_{p} =
k_{diss} \xi - \gamma \, \tau   
\left( k_{rec} + k_{diss} \mathrm{e}^{-t/\tau}\right)\, p \; n
+\ I,
\end{array}
\right.
\end{equation}
where 
\begin{equation}\label{eq:memory_term}
I :=
\gamma \, k_{diss}\ \displaystyle \int_{0}^{t} 
\left[p(\vect{x},s) n(\vect{x},s)-p(\vect{x},t) n(\vect{x},t) \right] 
\mathrm{e}^{-(t-s)/\tau} \ ds.
\end{equation}
Although $I$  is no more a convolution integral, it has the interesting 
property of vanishing both at $t=0$ and $t=+\infty$,
from which we expect, at least formally, that replacing the integral $I$ by a suitable approximation, say $\widetilde I$,
should not have a significant impact on the model behaviour as long as it preserves the asymptotics of $I$.
Our choice is to use a 
 trapezoidal quadrature rule, yielding
\begin{align}\label{eq:quadrule}
I \simeq
\widetilde{I} 
= \gamma \, k_{diss}\ \Frac{t}{2}\mathrm{e}^{-t/\tau}\left[
p(\vect{x},0) n(\vect{x},0)- p(\vect{x},t) n(\vect{x},t) \right].
\end{align}
It is easy to see that $\widetilde I$ vanishes both at $t=0$ and $t=+\infty$; moreover, the approximate formula~\eqref{eq:quadrule} as the advantage of \emph{lumping} the non-locality of $I$ into a quadratic term that has the same form as the generation/recombination rates already present in the right-hand side of~\eqref{eq:ncont_no_x_2} 
The resulting reduced model reads:
\begin{equation}\label{eq:model-transient}
\left\{
\begin{array}{lll}
-\Div(\varepsilon \nabu \varphi) & = & q (p - n) \\[2mm]
\Frac{\partial n}{\partial t}-\Div \vect{J}_n  & = & 
\widetilde G_{n} 
- \widetilde R_n
n \\[3mm]
\Frac{\partial p}{\partial t} 
-\Div \vect{J}_p & = & \widetilde G_{p} 
- \widetilde R_p 
p, \\[2mm]
\end{array}
\right. \qquad \mbox {in }\Omega_{T}
\end{equation}
where the modified generation/recombination mechanisms are
defined as
\begin{equation}\label{eq:eq_sys_const_rel}
\left\{
\begin{array}{ll}
\widetilde G_{n} 
& = \widetilde G_{p}
= 
k_{diss} \xi(\vect{x},t) + \gamma \; k_{diss} \Frac{t}{2} 
\mathrm{e}^{-t/\tau} p(\vect{x},0) n(\vect{x},0) \\[8mm]
\widetilde R_{n}
n  & = \widetilde R_{p}
p =
\gamma \left[ \tau (k_{rec}+k_{diss} \mathrm{e}^{-t/\tau}) 
+ k_{diss} \Frac{t}{2} \mathrm{e}^{-t/\tau}\right]  p(\vect{x},t) n(\vect{x},t) .
\end{array}
\right.
\end{equation}
Having derived a new, simplified model, 
it is natural to ask to which extent the novel formulation
is capable to describe correctly the main features of
the performance of an OSC. With this aim, we first investigate
the quality of the approximation provided by $\widetilde{I}$;
the quadrature error 
associated with the use of the
trapezoidal rule in~\eqref{eq:quadrule} is given by the
following relation~\cite{QSS2007}
\begin{equation}\label{eq:quadr_error}
E(t) = -\Frac{t^3}{12} e^{-(t-\zeta)/\tau}
\left( \lambda^{\prime \prime} (\zeta) +
\Frac{2}{\tau} 
\lambda^{\prime} (\zeta) + \Frac{1}{\tau^2} 
(\lambda(\zeta) - \lambda(t)) \right)
\end{equation}
where $\zeta \in (0,t)$ and $\lambda(s):=p(\cdot, s)\ n(\cdot, s)$.
Eq.~\eqref{eq:quadr_error}
shows that $E(t)$ becomes negligible as $t \rightarrow 0$ 
or $t \rightarrow +\infty$, as expected, meaning that 
the predicted (computed) stationary current is independent of 
the use of~\eqref{eq:quadrule} or the exact 
expression~\eqref{eq:ncont_no_x_2}, as 
numerically verified in Sect.~\ref{sec:num_constant_coeff}. 
However, for a finite value of time $t$, 
the discrepancy between the exact convolution term
and its approximation may be non-negligible. 
A reasonable estimate of the error would require a 
knowledge on the temporal behavior of the photogenerated
carrier densities $n$ and $p$ as a function of time.
This knowledge not being available, we can still
gain some information on the quadrature error by
an analogy with the approximation of the 
recombination/generation term that is usually 
carried out in the study of currents in a $p-n$ junction 
in the inorganic case (see~\cite{Sze2006}). This analogy
suggests that the value of $E(t)$ during the photocurrent 
transient (i.e., for $t$ sufficiently far from $0$
but also sufficiently far from stationary conditions)
might become significant if the OSC is operating under
high injection conditions, or, equivalently, 
high current level conditions. Again, this latter statement 
is numerically verified in Sect.~\ref{sec:num_constant_coeff}.

\begin{theorem}[Existence of a solution in the transient regime]
\label{th:transient}
Let assumptions (H1)--(H4) be satisfied, and 
the initial data $\vect{U}:=(n_0, p_0)$, $X_0$ and the function 
$\Psi$ be such that $\vect{U} \in (H^1(\Omega_T) \cap L^{\infty}(\Omega_T))^{2}$,
with $\vect{U} > \vect{0}$, $X_0 \in L^{\infty}(\Omega)$ 
with $X_0 \geq 0$, and $\Psi \in H^1(\Omega_T) \cap L^{\infty}(\Omega_T)$. 
Then, setting $\vect{u}: = (n, p)$, 
system~\eqref{eq:model-transient}-\eqref{eq:eq_sys_const_rel},
supplied with the initial/boundary conditions~\eqref{eq:electronboundaryconditions}--\eqref{eq:potentialboundaryconditions}, admits a weak solution 
$(\varphi, \vect{u})$ such that:
\begin{enumerate}
\item $\vect{u} > \vect{0}$ a.e. in $\Omega_T$;
\item $\vect{u}(\vect{x},0) = \vect{U}(\vect{x},0)$ and
$\vect{u}-\vect{U} \in L^2\left(0,T; H_{0}\right)^{2}$;
\item $\vect{u} \in \left( C(0,T; L^2(\Omega)) \cap L^{\infty}(\Omega_T) \right)^{2}$;
\item $\displaystyle 
\frac{\partial \vect{u}}{\partial t} \in L^2(0,T; {H}^{\prime}_{0})^{2}$;
\item $\varphi - \Psi \in L^2(0,T; H_0)$ with 
$\varphi \in L^{\infty}(\Omega_T)$,
\end{enumerate}
where $H_0 := \left\{ v \in H^{1}(\Omega) \; : \; v|_{\Gamma_D} = 0 \right\}$ and $H_{0}^{\prime}$ is its dual.\\
Moreover, using~\eqref{eq:x_di_t} and the 
regularity of $X_0$, $n$ and $p$, we have that 
$$
X, \Frac{\partial X}{\partial t} \in C(0,T; L^2(\Omega)) 
\cap L^{\infty}(\Omega_T)
$$ 
with $X(\vect{x},t) >0$ for all $t >0$ and for a.e. $\vect{x} \in \Omega$.
\end{theorem}

The proof of Theorem~\ref{th:transient} consists of verifying
that all of the assumptions (Ei)--(Eiv) of~\cite{Wrzosek1994}, p. 296
are satisfied.
It is immediate to see that the 
functions $\widetilde{R}_{\eta}$ are positive for $p>0$ and $n>0$ and 
satisfy locally Lipschitz conditions, with a Lipschitz constant 
which is uniform in time and equal to $2 \gamma$.
As a matter of fact, 
for any $n^{\prime},p^{\prime}$, $n^{\prime \prime},p^{\prime \prime}$
and for any $\vect{x}$ and $t$, we have
$$
\begin{array}{lll}
|\widetilde{R}_n(\vect{x},t,n^{\prime},p^{\prime})-
\widetilde{R}_n(\vect{x},t,n^{\prime \prime},p^{\prime \prime})|
& \leq &
\gamma \; \left( \tau (k_{rec}+k_{diss}) + \tau k_{diss} \right)
|p^{\prime} - p^{\prime \prime}| \\[2mm]
& \leq & 2 \gamma |p^{\prime} - p^{\prime \prime}|,
\end{array}
$$
and the same estimate holds for $\widetilde{R}_p$ provided to
exchange $|p^{\prime} - p^{\prime \prime}|$ with
$|n^{\prime} - n^{\prime \prime}|$.
Moreover, (H2) and (H3) ensure that 
\eqref{eq:model-transient}$_{2,3}$ are uniformly elliptic
with uniformly bounded convective velocities.
Then, by applying Theorem 2 of~\cite{Wrzosek1994}, we conclude
that Theorem~\ref{th:transient} (of the present article!) holds.

\section{Numerical Discretization}\label{sec:numericalmethod}
In this section, we illustrate the numerical techniques for the simulation 
of the full model~\eqref{eq:modelequations}--\eqref{eq:boundaryconditions}, as the
same approach can be used, with slight modifications, to treat
the reduced approximate 
model~\eqref{eq:model-transient}--\eqref{eq:eq_sys_const_rel}. 
In designing the algorithm presented here, our aim is twofold: on the one 
hand, it seems natural to try to adapt methods that are known to 
work efficiently and reliably for transient simulation of 
inorganic semiconductor devices (see, {\it e.g.}, 
\cite{Sel84} Chapt.~6, Sect.~4); on the other hand, as the emphasis 
of the present paper is on accurately estimating photocurrent 
transient times, it is necessary to apply advanced time-step control 
techniques~\cite{ascher1998cmo,Hairer96b}.
To this end, our chosen approach is based on Rothe's method 
(also known as method of horizontal lines) which consists 
of three main steps: first, the time dependent problem 
is transformed into a sequence of stationary differential problems by 
approximating the time derivatives by a suitable 
difference formula; then, the resulting nonlinear 
problems are linearized by an appropriate functional 
iteration scheme; and, finally, the linear 
differential problems obtained are solved numerically using a 
Galerkin--Finite Element Method (G--FEM) for the spatial discretization. 
Sects.~\ref{sec:timedisc},~\ref{sec:lin} and~\ref{seq:spacedisc} 
below discuss in more detail each of these steps; it is worth noting
that, with minor modifications, the linearization techniques 
of Sect.~\ref{sec:lin} can also be applied to treat the 
stationary model~\eqref{eq:model-stationary}.

\subsection{Time Discretization}\label{sec:timedisc}
To transform the time dependent 
problem~\eqref{eq:modelequations}--\eqref{eq:boundaryconditions} 
into a sequence of stationary problems, we replace
the partial time derivative with a suitable finite difference 
approximation, specifically, the Backward Differencing Formulas (BDF) 
of order $m \le 5$ (see, {\it e.g.},~\cite{ascher1998cmo}, Sect.~10.1.2).
To describe the resulting stationary problem, 
let $ 0 = t_0 <  \ldots < t_{K-1} < t_K < T $
be a strictly increasing, not necessarily uniformly 
spaced, finite sequence of time levels and assume the quantities
$u_1=n$, $u_2=p$, $X$ and $\varphi$ to be known functions of 
$\vect{x}$ for every $t_k$, $k=0\ldots K-1$. Then we 
obtain:
\begin{equation}\label{eq:td_general_model}
\left\{
\begin{array}{lll}
-\Div(\varepsilon \nabu \varphi_K) + q\ (n_K - p_K) & = & 0 \\[2mm]
\Sum{k=0}{m}  \theta_k n_{K-k} - \Div \vect{J}_n(n_K; 
\nabu \varphi_K) - U_K & = & 0\\[4mm]
\Sum{k=0}{m}  \theta_k p_{K-k} - \Div \vect{J}_p(p_K; 
\nabu \varphi_K) - U_K & = & 0\\[4mm]
\Sum{k=0}{m}  \theta_k X_{K-k} - W_K & = & 0,
\end{array}
\right.
\end{equation}
where $f_k = f (\vect{x}, t_k)$ for any generic function 
$f = f(\vect{x}, t)$, and 
\begin{eqnarray*}
U_K && := U (\nabu \varphi_K, n_K, p_K, X_K, t_K)\\ 
&& = G_n(\nabu \varphi_K, n_K, p_K, X_K, t_K) - 
R_n(\nabu \varphi_K, n_K, p_K, X_K, t_K)\ n_K\\
&& = G_p(\nabu \varphi_K, n_K, p_K, X_K, t_K) - 
R_p(\nabu \varphi_K, n_K, p_K, X_K, t_K)\ p_K,\\[5mm]
W_K && := W (\nabu \varphi_K, n_K, p_K, X_K, t_K)\\
&& = g (\nabu \varphi_K, n_K, p_K, X_K, t_K) - 
 r(\nabu \varphi_K, n_K, p_K, X_K, t_K).
\end{eqnarray*}
System~\eqref{eq:td_general_model}, 
together with the constitutive relations for the 
fluxes given in~\eqref{eq:driftdiffusion} and the 
set of boundary conditions~\eqref{eq:boundaryconditions},
constitutes a system of nonlinear elliptic differential 
equations \eqref{eq:continuity} coupled to an 
algebraic constraint equation \eqref{eq:poisson}.
In our implementation, the selection of the next time level $t_{K}$ 
and of the formula's order $m$, as well as the computation of the 
corresponding coefficients $\theta_{k}, k=0,\ldots,m$, is performed 
adaptively to minimize the time discretization error while 
minimizing the total number of time steps via the 
DAE solver software library DASPK~\cite{brown1992dds,vankeken1995dnh}.
Notice that, if $m=1$, we have  
$\theta_{0} = - \theta_{1} = 
\frac{1}{t_{K}-t_{K-1}}$, $\theta_{k} = 0,\; k > 1$, 
and the temporal semi-discretization of 
system~\eqref{eq:modelequations}--\eqref{eq:boundaryconditions} coincides with
the Backward Euler method. 

\subsection{Linearization}\label{sec:lin}
To ease the notation, throughout this section
the subscripts denoting the current time level will be dropped.
Let $\vect{y}:=[\varphi, \; n, \; p, \; X]^T$ denote the
vector of dependent variables and let $\vect{0}$ denote the
null vector in $\mathbb{R}^4$. Then, the nonlinear system~\eqref{eq:td_general_model}
can be written in compact form as 
\begin{equation}\label{eq:tdabs_general_model}
\vect{F}(\vect{y}) = \vect{0}, \qquad
\mbox{with } \qquad
\vect{F}(\vect{y}) =
\left\{
\begin{array}{l}
f _\varphi (\varphi, n, p)\\
f _n (\varphi, n, p, X)\\
f _p (\varphi, n, p, X)\\
f _X (\varphi, n, p, X)\\
\end{array}
\right\}.
\end{equation}
The adopted functional iteration technique for the linearization
and successive solution of problem~\eqref{eq:td_general_model} 
is the Newton-Raphson method. One step of this scheme 
can be written as
\begin{equation}\label{eq:tdlin_general_model}
\left[\begin{array}{cccc}
\partial_\varphi (f_\varphi) &
	\partial_n (f_\varphi) &
	\partial_p (f_\varphi)&
	0\\
\partial_\varphi (f_n) &
	\partial_n (f_n)&
	\partial_p (f_n)&
	\partial_X (f_n)\\
\partial_\varphi (f_p)&
	\partial_n (f_p)&
	\partial_p (f_p)&
	\partial_X (f_p)\\
\partial_\varphi (f_X)&
	\partial_n (f_X)&
	\partial_p (f_X)&
	\partial_X (f_X)
\end{array}\right]_{(\varphi, n, p, X)}
\left[
\begin{array}{l}
\Delta \varphi\\
\Delta n\\
\Delta p\\
\Delta X\\
\end{array}
\right] =
\left[\begin{array}{l}
-f _\varphi (\varphi, n, p) \\
-f _n (\varphi, n, p, X) \\
-f _p (\varphi, n, p, X) \\
-f _X (\varphi, n, p, X) \\
\end{array}\right]
\end{equation}
where $\partial_a (f)$ denotes the Frech\'et derivative of 
the nonlinear operator $f$ with respect to the function $a$.
More concisely, we can express~\eqref{eq:tdlin_general_model} 
in matrix form as
$$
\vect{J} (\vect{y}) \; \Delta \vect{y} = -
\vect{F}(\vect{y}),
$$
where $\vect{J}$ is the Jacobian matrix and 
$\Delta \vect{y}:=[\Delta \varphi, \; \Delta n, \; 
\Delta p, \; \Delta X]^T$ is the unknown increment vector.
The exact computation of all the derivatives in the Jacobian on the left 
hand side in~\eqref{eq:tdlin_general_model} can become quite complicated 
if the full model for all the coefficients (most notably the electric field 
dependence of $k_{diss}$, $\mu_n$ and $\mu_p$) is taken into account.
Moreover, this would require cumbersome modifications to the solver code 
whenever a new coefficient model is to be implemented. 
One alternative could be to employ a staggered solution algorithm, 
often referred to as Gummel-type approach 
in the semiconductor simulation context~\cite{gummel64,CdFJWJRS2009}. 
The decoupled approach is well known to be more robust as compared
to the fully coupled Newton approach~\eqref{eq:tdlin_general_model} with 
respect to the choice the initial guess and also less memory consuming.
As in this particular study we can rely on the knowledge of the system 
variables at previous time levels to construct a reasonable initial 
guess and as we are dealing with an intrinsically one-dimensional problem
(see Sect.~\ref{sec:numericalresults}), 
memory occupation is not likely to be a stringent 
constraint, so that we adopt a quasi-Newton method 
where, rather than the exact Jacobian $\vect{J}(\vect{y})$,
we use an approximation $\widetilde {\vect{J}}(\vect{y})$ 
in which the dependence of the mobilities, of the diffusion coefficients 
and of the dissociation coefficient on the electric field is neglected. 
This approach has the further advantage of facilitating the use of a standard 
software library like DASPK for advancing in time.

\subsection{Spatial Discretization and Balancing of 
the Linear System}\label{seq:spacedisc}
Once the linearization described in the previous 
section is applied, the resulting linear system of PDEs 
is numerically approximated by means of a suitable G--FEM. 
Precisely, to avoid instabilities and 
spurious oscillations that may arise when the drift 
terms become dominant, we employ an exponential fitting finite element 
discretization~\cite{Bank98,ricsac98,XuZikatanov1999,LazarovZikatanov2005}.
This formulation provides a natural multidimensional
extension of the classical Scharfetter-Gummel difference 
scheme~\cite{ScharfetterGummel1969,BrezziCMAME89} 
and ensures, when applied to a carrier continuity equation 
in the DD model, that the computed carrier concentration 
is strictly positive under the condition that the triangulation 
of the domain $\Omega$ is of Delaunay type.
It is important to notice that, 
when implementing on the computer the above described
procedure, the different physical nature of the unknowns of the system 
and their wide range of variation may lead to badly scaled and 
therefore ill-conditioned linear algebraic problems, which in turn 
can negatively affect the accuracy and efficiency of the algorithm.
To work around this issue, we introduce two sets of scaling coefficients,
denoted $\{ \sigma_{\varphi},\ \sigma_n,\ \sigma_p,\ \sigma_X\}$ and 
$\{\bar \varphi,\ \bar n,\ \bar p,\ \bar X\}$, and restate 
problem~\eqref{eq:tdabs_general_model} as 
\begin{equation}\label{eq:tdabs_general_model_sc}
\left\{\begin{array}{l}
\Frac{1}{\sigma_{\varphi}} f _\varphi 
(\bar \varphi \hat \varphi, \bar n \hat n, \bar p \hat p) = 0\\[3mm]
\Frac{1}{\sigma_{n}} f _n (\bar \varphi \hat \varphi, \bar n \hat n, 
\bar p \hat p, \bar X \hat X) = 0\\[3mm]
\Frac{1}{\sigma_{p}} f _p (\bar \varphi \hat \varphi, \bar n \hat n, 
\bar p \hat p, \bar X \hat X) = 0\\[3mm]
\Frac{1}{\sigma_{X}} f _X (\bar \varphi \hat \varphi, \bar n \hat n, 
\bar p \hat p, \bar X \hat X) = 0,
\end{array}\right.
\end{equation}
where $\hat \varphi:= \varphi/\bar \varphi$, 
$\hat n := n/\bar n$, $\hat p := p/\bar p$ and 
$\hat X := X/\bar X$. Solving~\eqref{eq:tdabs_general_model_sc} for the 
scaled dependent variables 
$[\hat \varphi, \; \hat n, \; \hat p, \; \hat X]^T$
corresponds to solving a system equivalent 
to~\eqref{eq:tdlin_general_model} where the rows 
of the Jacobian $\vect{J}$ and of the residual $\vect{F}$ 
are multiplied by the factors 
$\{ 1/\sigma_{\varphi},\ 1/\sigma_n,\ 1/\sigma_p,\ 1/\sigma_X \}$ 
while the columns of $\vect{J}$ are multiplied by 
the factors $\{\bar \varphi,\ \bar n,\ \bar p,\ \bar X\}$.
Computational experience reveals that a proper choice of 
the scaling coefficients might have a strong impact 
on the performance of the algorithm. 
For example, to obtain the results of Fig.~\ref{fig:numres_fig6} 
a suitable choice was found to be that of setting
$\sigma_{\varphi} = 1,\ \sigma_n=\sigma_p,=10^{3},\ 
\sigma_X=10^{2}$ and $\bar {\varphi} = 1,\ 
\bar n=\bar p,=10^{22},\ \bar X=10^{19}$ 
while values differing by more than one order 
of magnitude from such choice were found to hinder 
the ability of the DAE solver to reach convergence.

\section{Numerical Results}\label{sec:numericalresults}

This section is devoted to presenting the results of numerical 
simulations carried out with the algorithms described in 
Sect.~\ref{sec:numericalmethod}. In particular, in 
Sect.~\ref{sec:num_varying_coeff} we discuss the simulation results 
for a realistic BHJ device focusing on the impact of 
the model parameter values on the turn-on transient time in different 
operation conditions. Sect.~\ref{sec:num_constant_coeff} is devoted 
to characterizing the region in the model parameter space where the 
approximate formula~\eqref{eq:quadrule} and the resulting reduced 
model~\eqref{eq:model-transient} are reliable. 
In both cases the considered device has a thickness $L_{OSC} = 70 nm$ 
and the contact materials are ITO and Al for the transparent 
and reflecting contact, respectively. As no external voltage is applied to the 
device, this results in a total voltage drop across the device $\Delta V = 0.5 V$. 
The relative permittivity constant is $\varepsilon_{r} = 4$ and the operating 
temperature is $300 K$. As the thickness of the device is much smaller with respect to the dimensions in the other directions (typically many orders of magnitude larger) and the donor/acceptor blend is considered to be uniform, the simulations presented here are performed in one spatial dimension, so that the computational domain is modeled as the segment $\Omega = [0, L_{OSC}]$ with the cathode at $x=0$ and the anode at $x=L_{OSC}$. Also, as the device length is quite small compared to the wavelength of visible light, it is reasonable to consider the photon absorption rate $G$ to be constant in $\Omega$ at any $t \in [0, T]$.

\subsection{Simulation of a realistic device}\label{sec:num_varying_coeff}

\begin{figure}[h]
\begin{center}
\subfigure[]{\includegraphics[width=.45\linewidth]{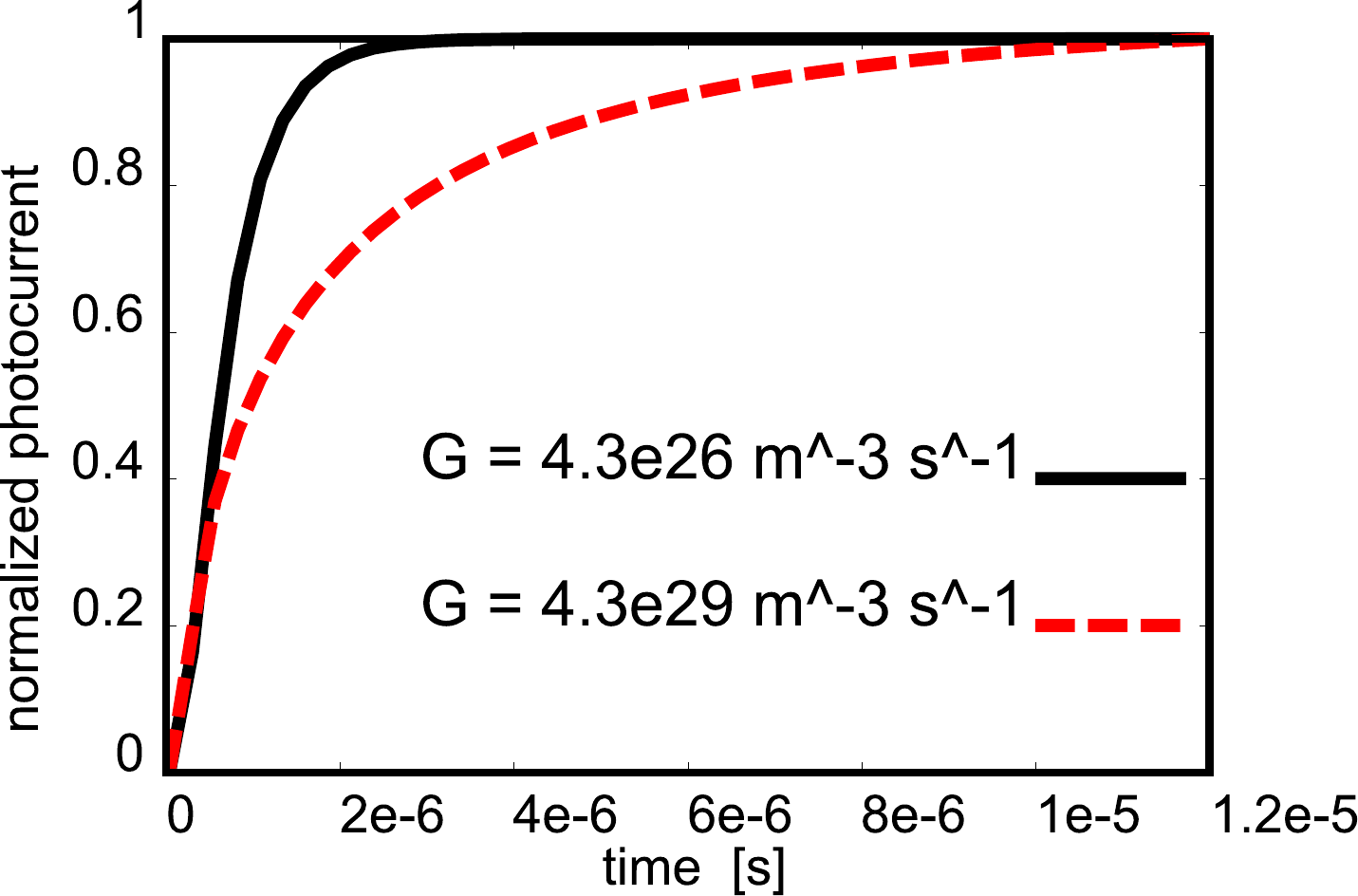}\label{fig:2a}}
\subfigure[]{\includegraphics[width=.45\linewidth]{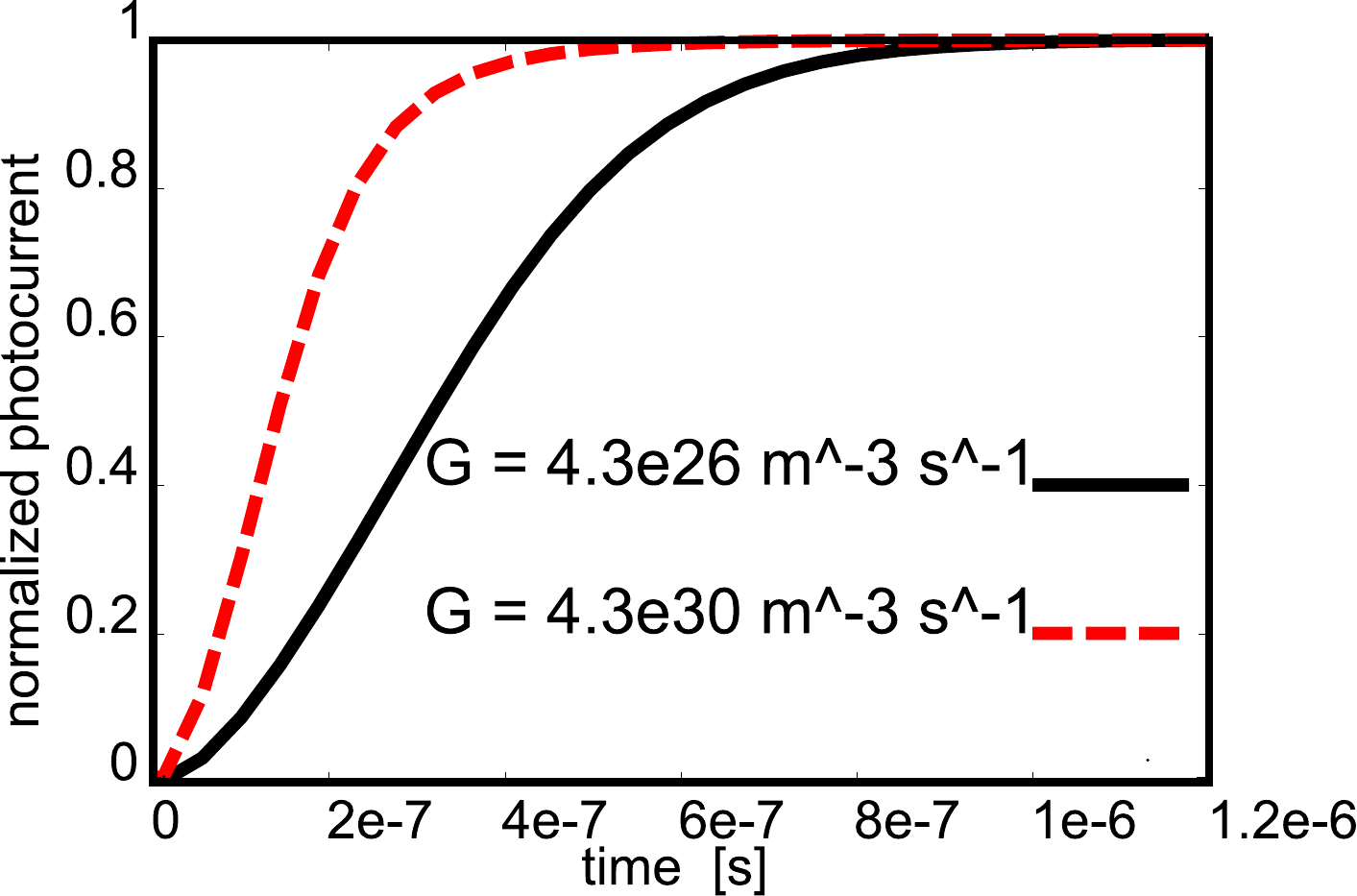}\label{fig:2b}}
\subfigure[]{\includegraphics[width=.45\linewidth]{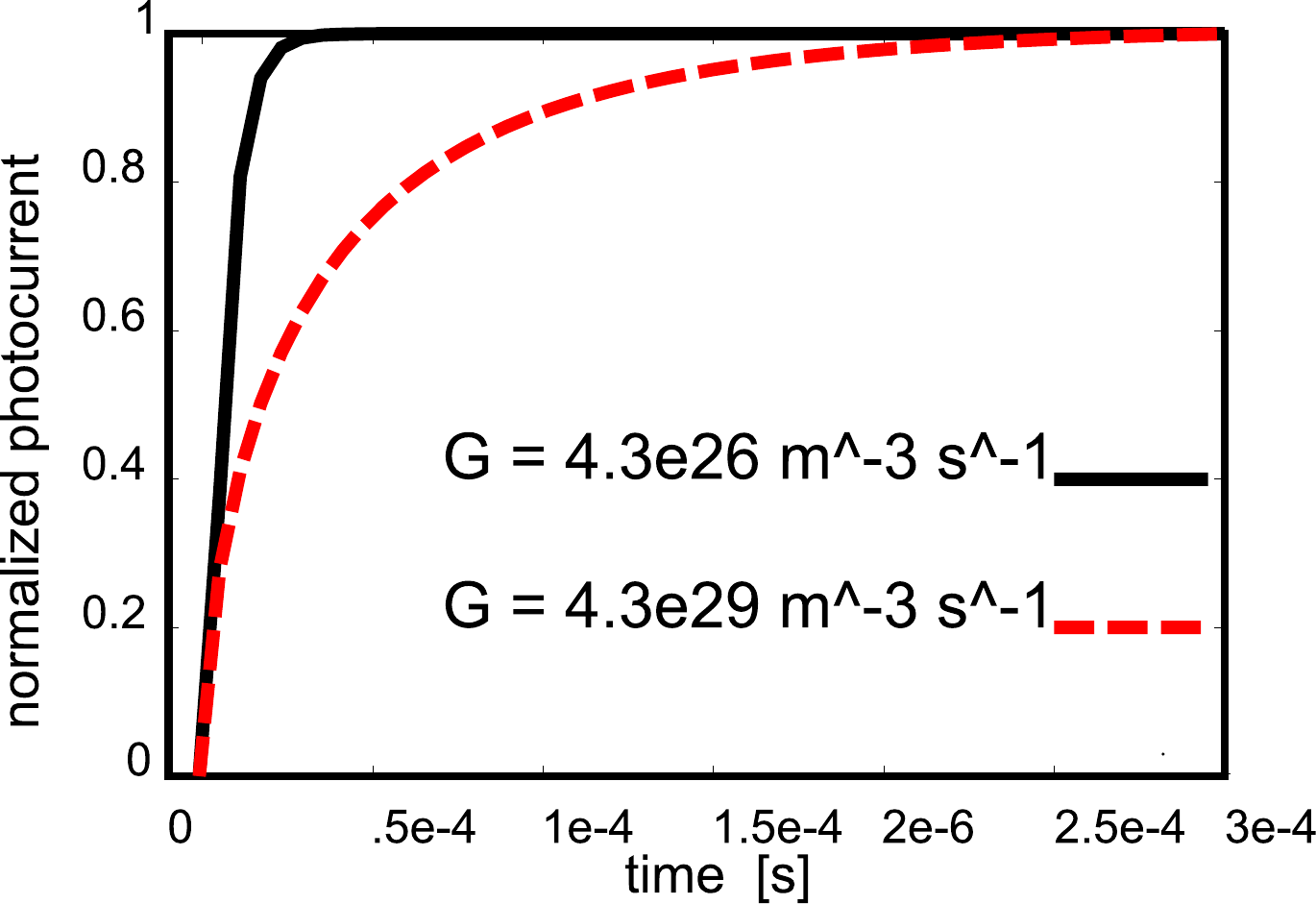}\label{fig:2c}}
\subfigure[]{\includegraphics[width=.45\linewidth]{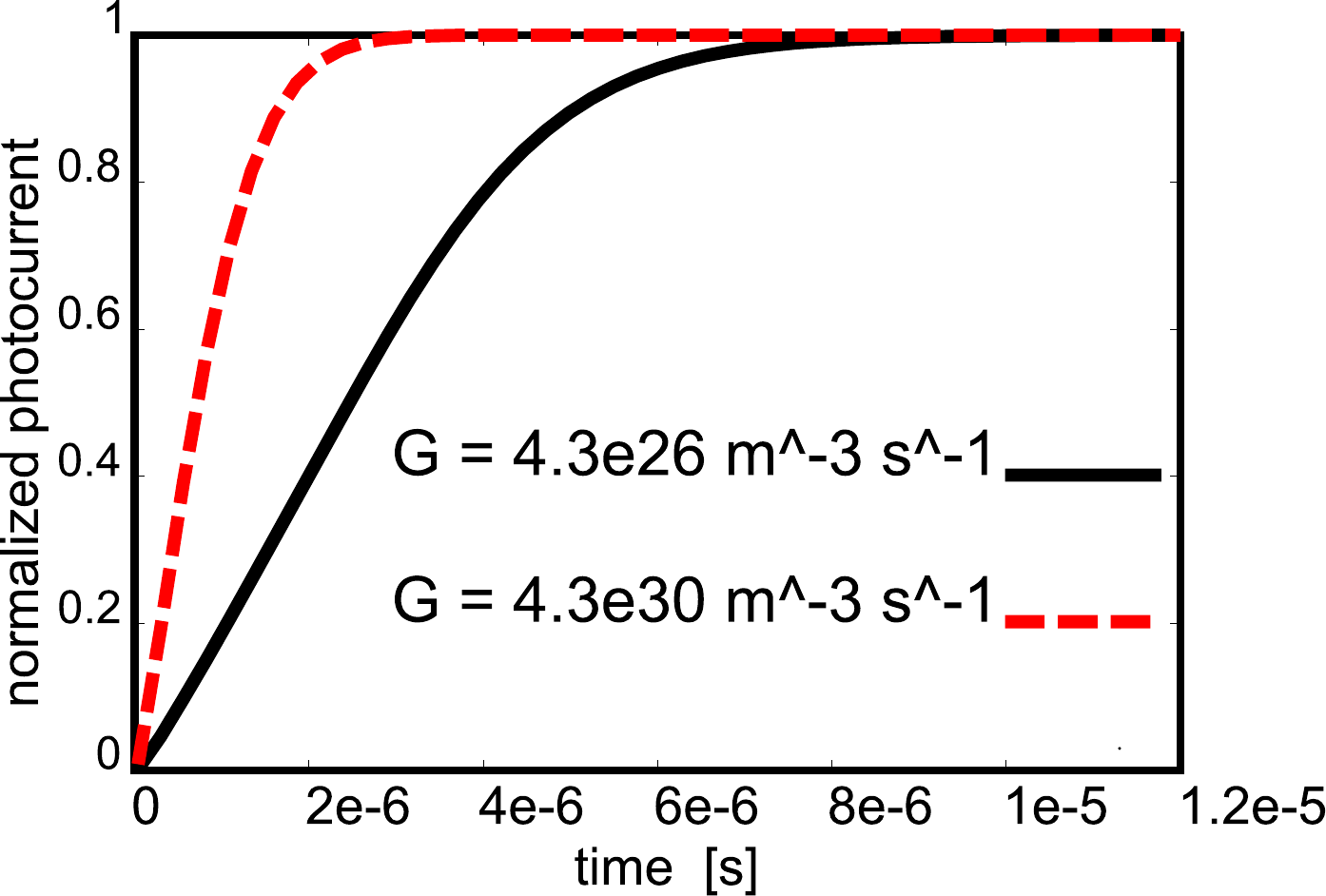}\label{fig:2d}}
\caption{Transient currents at low and high intensities 
with different mobilities and exciton recombination rate coefficients. 
For (a) and (b) the mobility was 
$2 \times 10^{-4} cm^{2} V^{-1} s^{-1}$ with 
geminate recombination rate constants $k_{rec} = 1 \times 10^{ 5} s^{-1}$ 
and $1 \times 10^{ 7} s^{-1}$ respectively. 
For (c) and (d) the mobility is 
$2 \times 10^{-5} cm^{2} V^{-1} s^{-1}$ with $k_{rec} = 1 \times 10^{ 4} s^{-1}$ 
and $k_{rec} = 1 \times 10^{ 6} s^{-1}$ respectively.  
}
\label{fig:numres_fig2}
\end{center}
\end{figure}

In this section we present simulation results of the realistic BHJ
device whose data are given in~\cite{Hwang2008}.
The computations try to reproduce the measurements that are commonly 
performed in research laboratories to characterize the device material
properties and are meant to show the ability of the model to capture
the complex dependence of the turn-on transient time on both the mobility
coefficients and the exciton dissociation/recombination dynamics, and the predominance 
of one or the other of such phenomena depending on the operation conditions, {\it i.e.}
on the intensity of the light to which the device is exposed. 
Throughout this section we use for the coefficients in the boundary conditions~\eqref{eq:boundaryconditions} 
the current injection model of~\cite{CampbellScott1999,Barker2003} corrected 
as in~\cite{Hwang2008,Hwang2009_private} to increase the carrier surface recombination rate, thus avoiding 
the occurrence of spurious charge build-up effects near the contacts.
The exciton dissociation coefficient $k_{diss}$ is considered to depend on the electric field according to an 
Onsager-like model given by the nonlinear formula presented in~\cite{Barker2003} with the initial separation of the geminate pair set to $a=1.5$nm, while the recombination $k_{rec}$ rate
is constant. The bimolecular recombination coefficient $\gamma$ depends on the carrier mobilities and on the material permittivity $\varepsilon$ as resulting from Langevin theory~\cite{Mihailetchi2004}, therefore, as we consider here 
the carrier mobilities to be constant, $\gamma$ is a constant as well.
Figure~\ref{fig:numres_fig2} shows the photocurrent evolution in response to an abrupt turn-on of a light source;
for each row in the figure the charge carrier mobilities are kept constant while the exciton recombination coefficient  
is varied whereas for each row in the figure the mobilities are fixed and the recombination coefficients vary.
By comparing Figs.~\ref{fig:2a} and~\ref{fig:2c} to Figs.~\ref{fig:2b} and~\ref{fig:2d} one can notice that the strong impact
of the recombination rate coefficient $k_{rec}$ on the transient duration in high illumination conditions (dashed lines) completely overshadows the effect of the carrier transport properties, while in low illumination conditions (solid lines) the importance of the effect of $k_{rec}$ is less apparent so that the transient time is more related to the value of the mobilities.

\begin{figure}[h]
\begin{center}
\subfigure[$k_{rec} = 10^{5} s^{-1}$]
{\includegraphics[width=.45\linewidth]{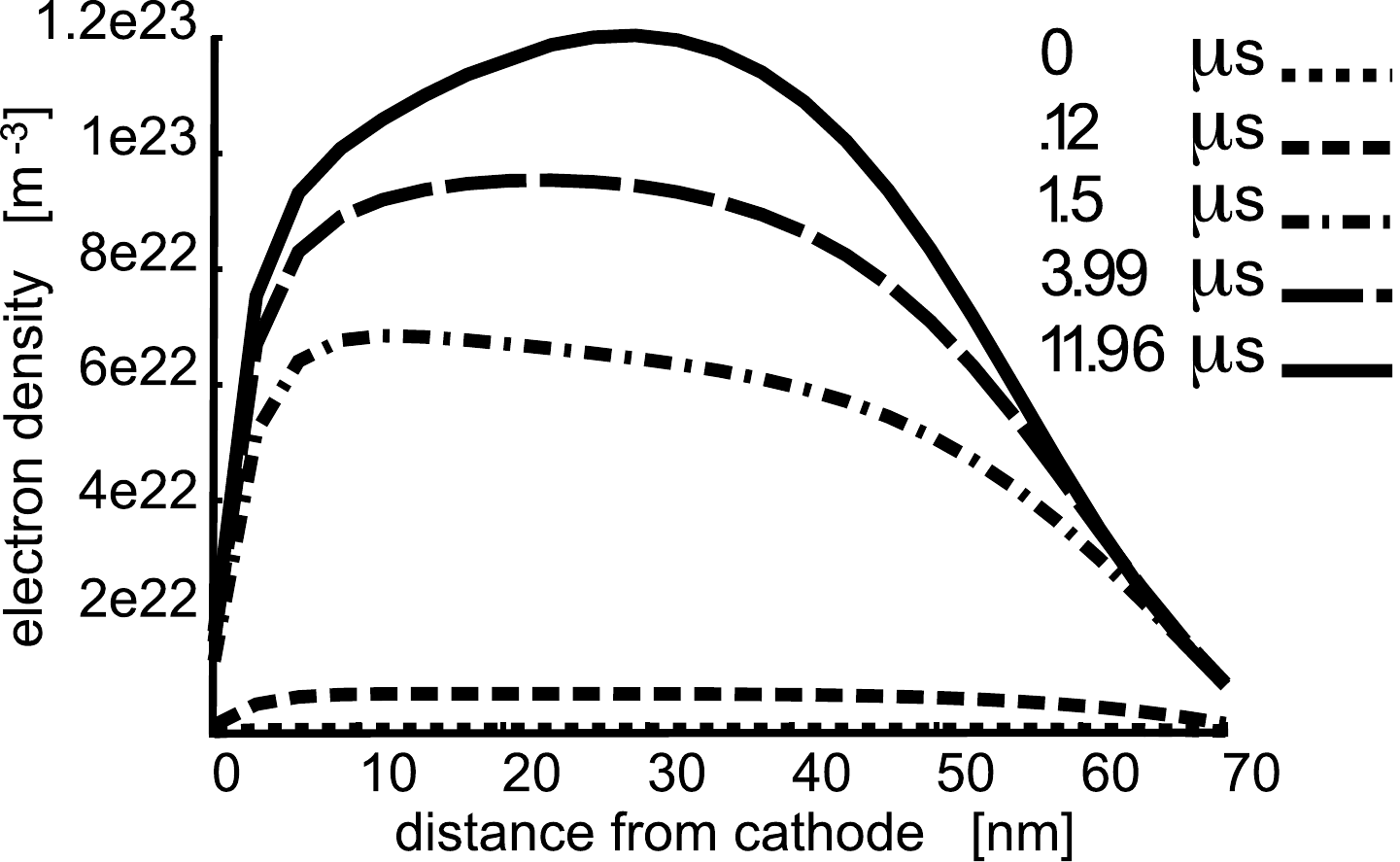}}
\subfigure[$k_{rec} = 10^{7} s^{-1}$]
{\includegraphics[width=.45\linewidth]{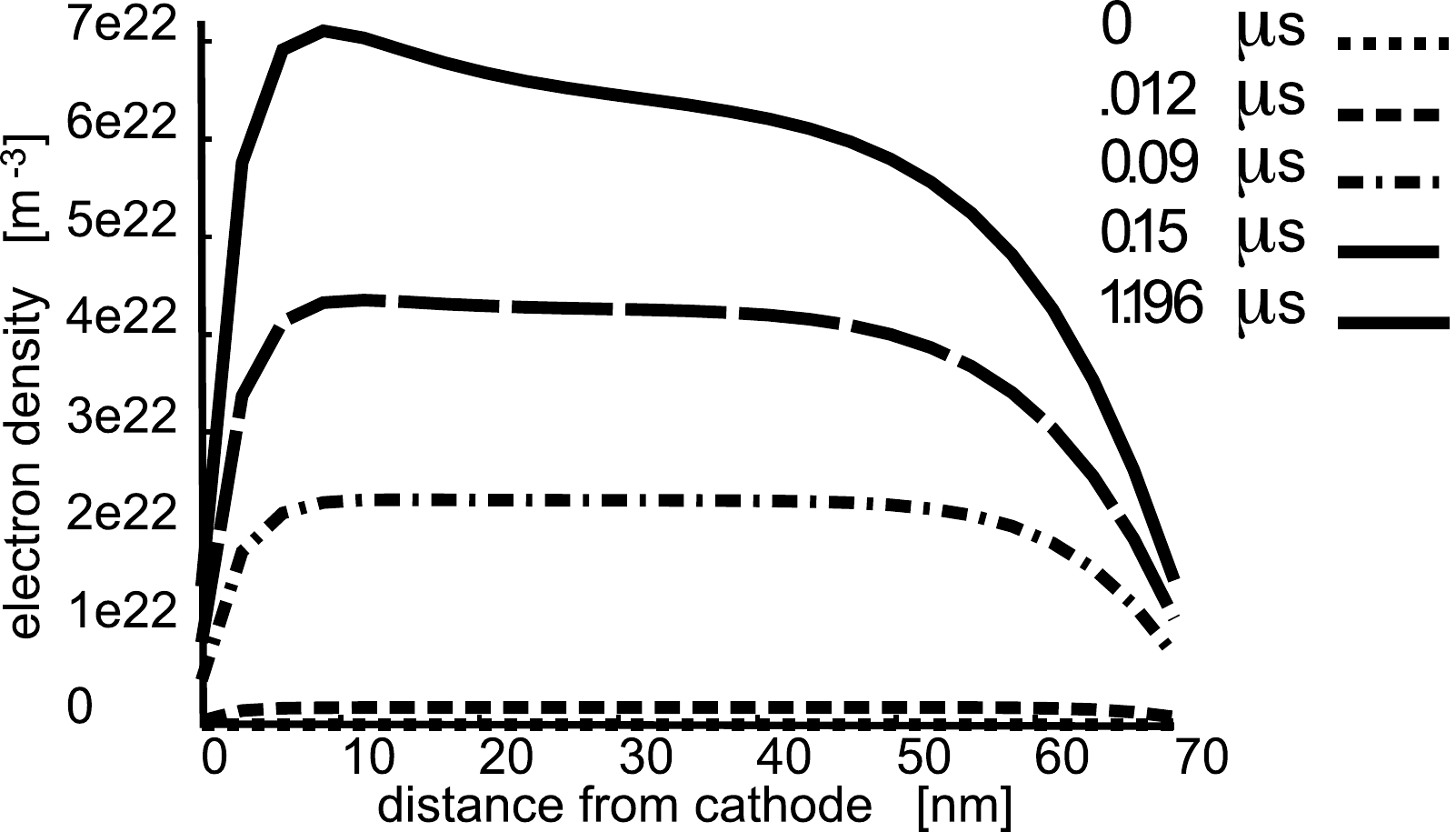}}
\caption{Time evolution of the electron distribution 
at high intensity with (a) high charge generation efficiency 
and (b) low charge generation efficiency.} 
\label{fig:numres_fig6}
\end{center}
\end{figure}

Figure~\ref{fig:numres_fig6} shows the time evolution of 
the electron density in the device under strong illumination 
conditions ($G = 4.3 \cdot 10^{30} \; m^{-3} s^{-1}$).
Hole density is not shown in the figures because, due to 
the choice of equal mobilities, it is the exact mirror 
image of the electron density. 
As previously mentioned, due to the absence of fixed charges (dopants)
within the bulk of the device the charge densities do not show the steep
interior layers that are the main peculiarity of inorganic semiconductor
models and lead to the main difficulties in the numerical simulation of such devices.
Also the steepness of the boundary layers is less extreme in the case of organic 
devices and is further mitigated by the inclusion of finite surface recombination speed 
in the boundary conditions.

\begin{figure}[h]
\begin{center}
\includegraphics[width=.45\linewidth]{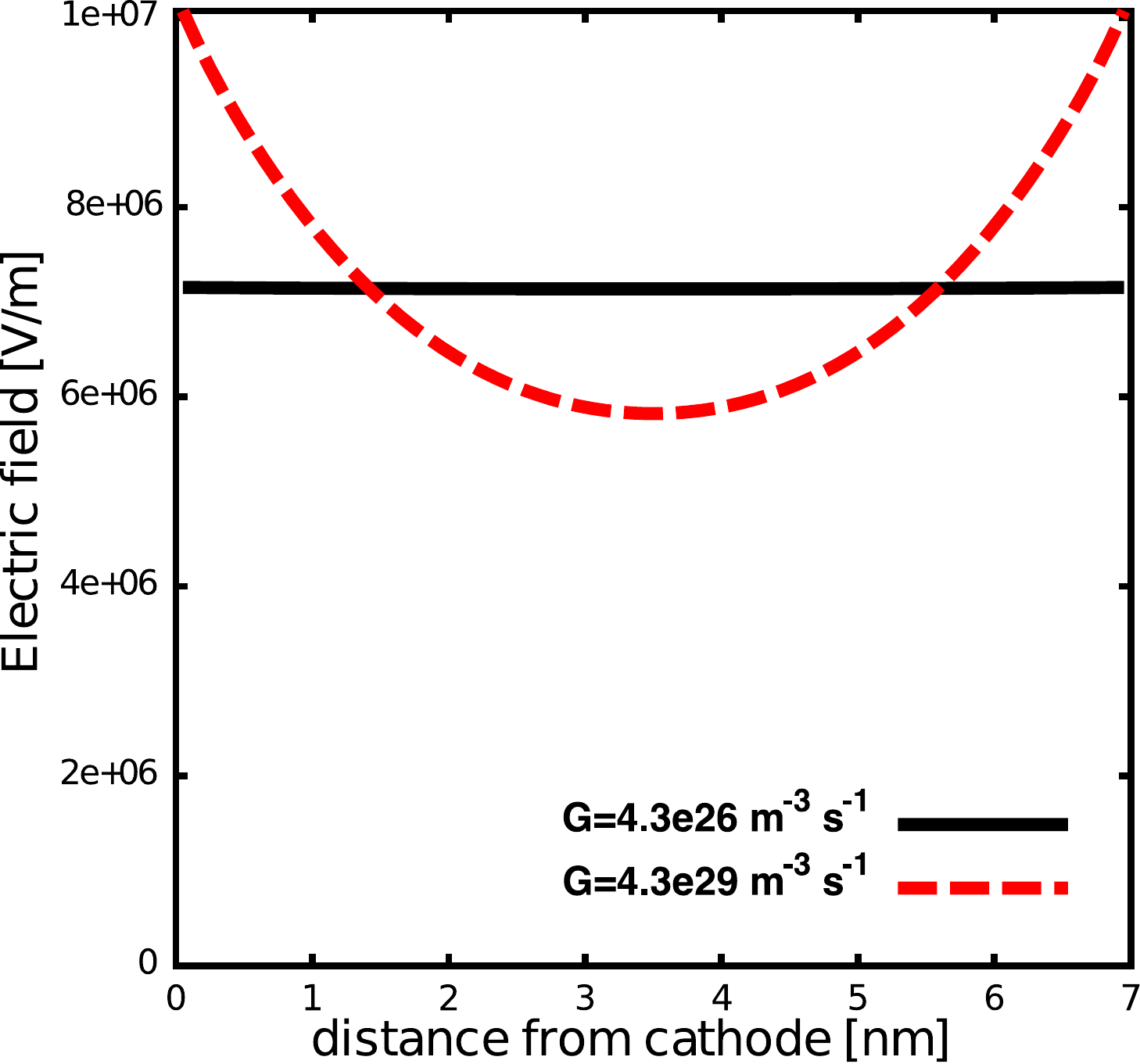}
\caption{Value of the computed electric field for a device with mobilities 
$\mu_{\eta} = 2 \times 10^{-4} cm^{2} V^{-1} s^{-1}$ and recombination 
rate constants $k_{rec} = 1 \times 10^{ 5} s^{-1}$.}
\label{fig:numres_fig_5_ab}
\end{center}
\end{figure}

The consistency of the results shown here with those of~\cite{Hwang2008}
is a strong indication of the robustness of the numerical algorithm of Sect.~\ref{sec:numericalmethod}.
Finally, Figure~\ref{fig:numres_fig_5_ab} shows the magnitude of the electric field 
along the device,
for low illumination (solid line) the electric field is practically constant
throughout the device while for high light intensity (dashed line) its deviation
around its mean value $\left<E\right> = \Delta V / L_{OSC}$ is about $30$\%.

\subsection{Validation of the simplified model}\label{sec:num_constant_coeff}

In this section we wish to estimate the impact of the approximation~\eqref{eq:quadrule}
on the simulation results for parameter values within a physically plausible range. 
To be consistent with assumptions (H1)-(H4) of Sect.~\ref{sec:analysis}, throughout the
present section we enforce that all model coefficients be constant by 
replacing the spatially varying electric field $E$ in the coefficient models
 by its mean value $\left<E\right> =-\Delta V /L_{OSC}$. Furthermore we  consider carrier recombination at the contacts be instantaneous, so that the boundary conditions~\eqref{eq:boundaryconditions} degenerate into simple Dirichlet type conditions. 
The plausibility of these assumptions
has been already addressed at the beginning of Sect.~\ref{sec:analysis} and in the discussion 
of the numerical results of Sect.~\ref{sec:num_varying_coeff}.
In all subsequent figures, the dashed line refers to
the solution computed with the full (3 carrier)
model~\eqref{eq:modelequations}-\eqref{eq:boundaryconditions}
while the solid line refers to the simplified approximate (2 carrier)
model~\eqref{eq:model-transient}-\eqref{eq:eq_sys_const_rel}. 

\begin{figure}[h]
\subfigure[Low mobilities: $\mu_{n}=\mu_{p} = 
2\times 10^{-9} m^{2} V^{-1} s^{-1}$]
{\includegraphics[width=.45\linewidth]{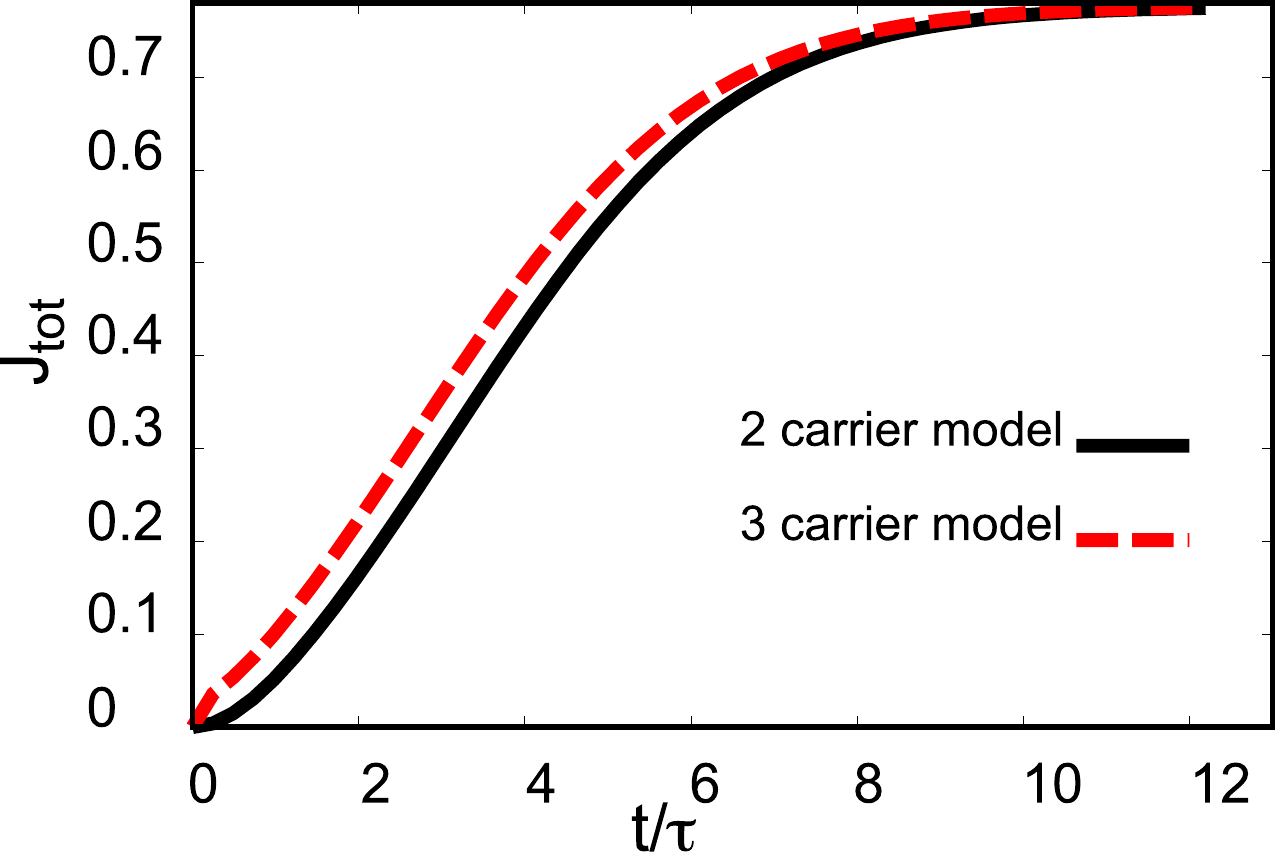}}
\subfigure[High mobilities: $\mu_{n}=\mu_{p} = 
2\times 10^{-8} m^{2} V^{-1} s^{-1}$]
{\includegraphics[width=.45\linewidth]{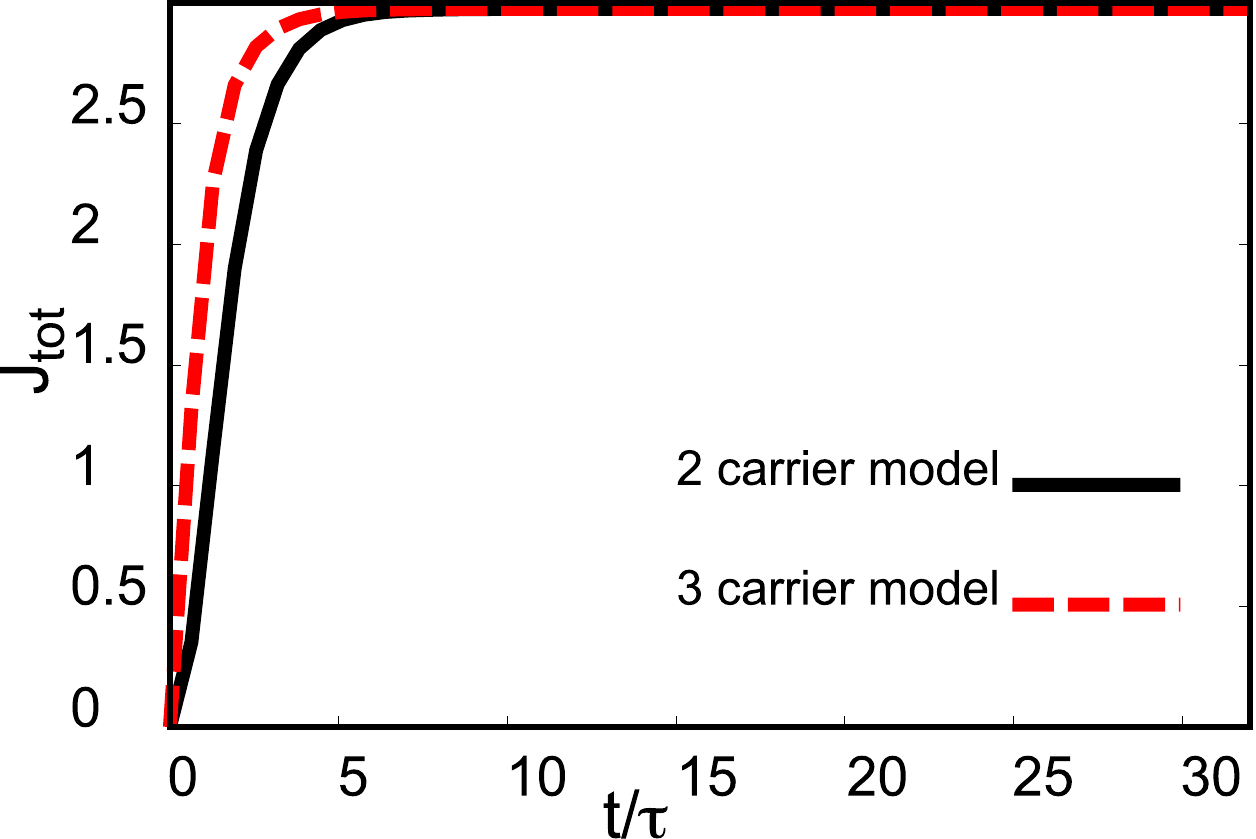}}
\caption{Photocurrent transient at low light intensity: 
effect of mobility on rise time.}
\label{fig:G-u}
\end{figure}

\begin{figure}[h]
\subfigure[$k_{diss}=4.4 \times 10^5 s^{-1}$]
{\includegraphics[width=.45\linewidth]{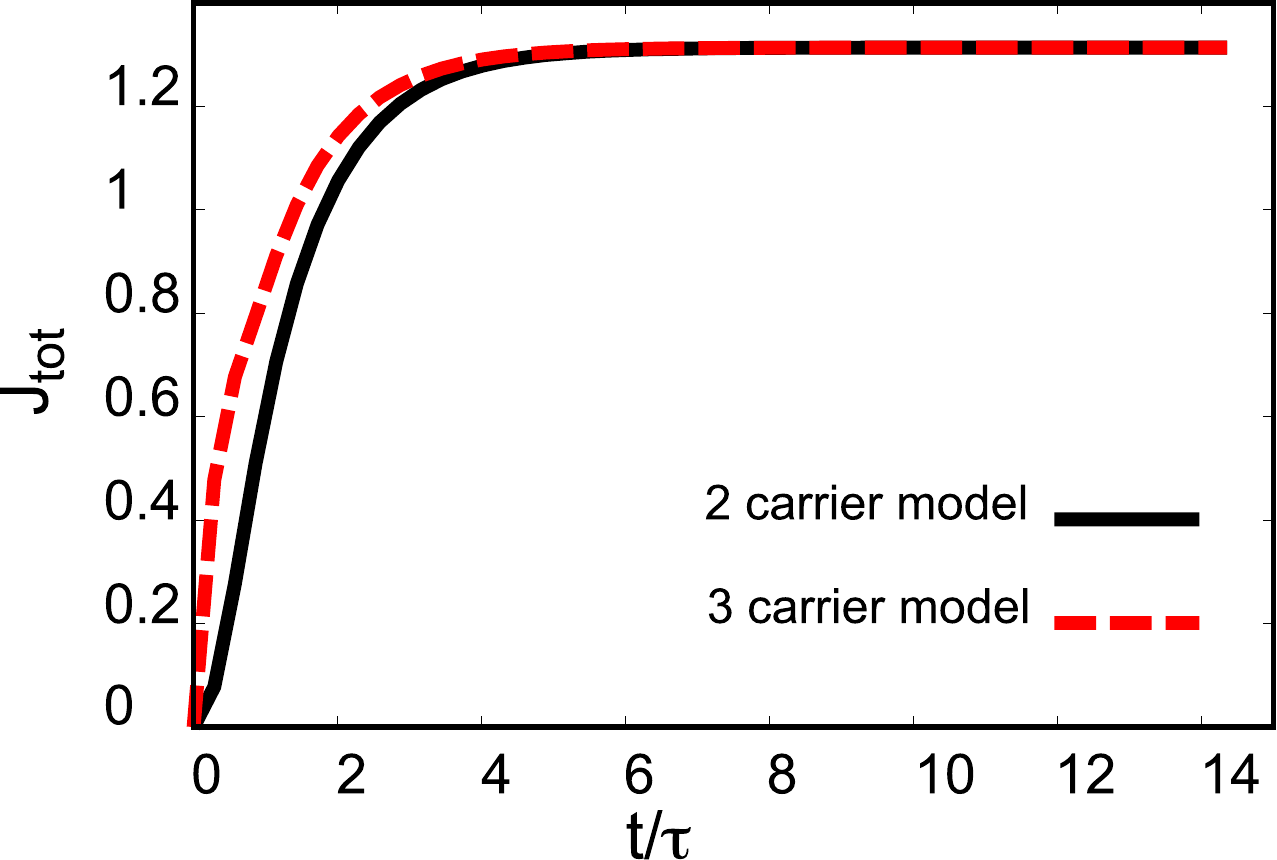}}
\subfigure[$k_{diss}=8 \times 10^6 s^{-1}$]
{\includegraphics[width=.45\linewidth]{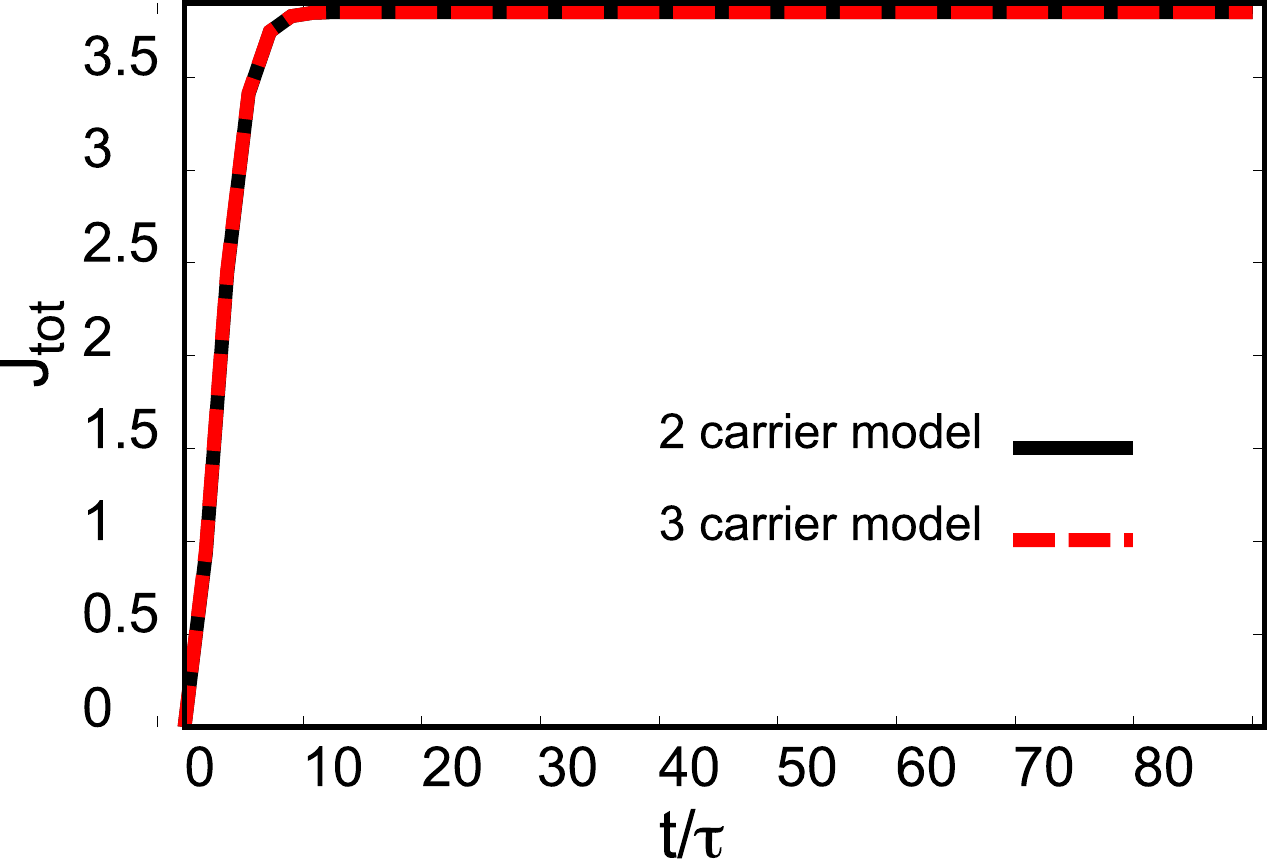}}
\caption{Photocurrent transient at low light intensity: effect 
of dissociation rate on rise time.}
\label{fig:G-kd}
\end{figure}

\begin{figure}[h]
\subfigure[$k_{rec}=10^5 s^{-1}$]
{\includegraphics[width=.45\linewidth]{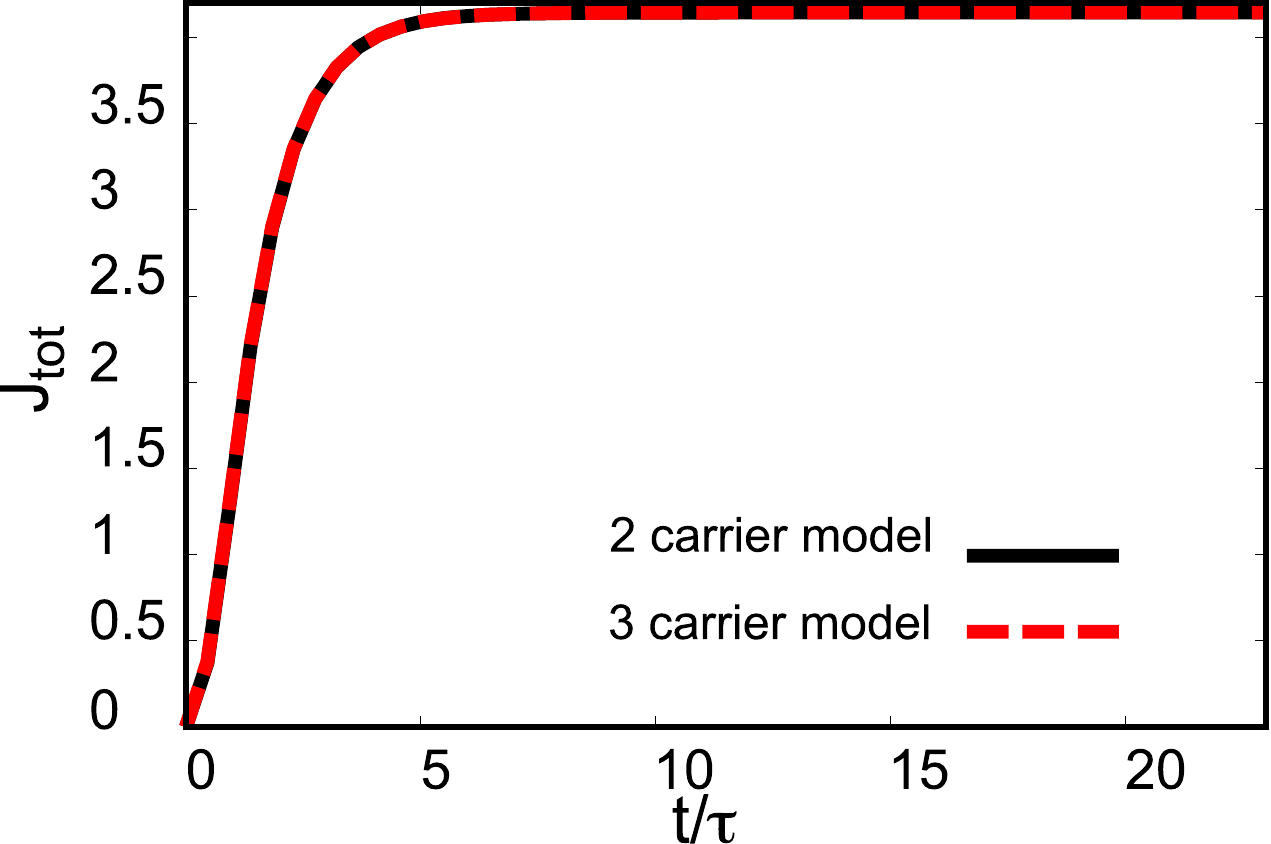}}
\subfigure[$k_{rec}=10^7 s^{-1}$]
{\includegraphics[width=.45\linewidth]{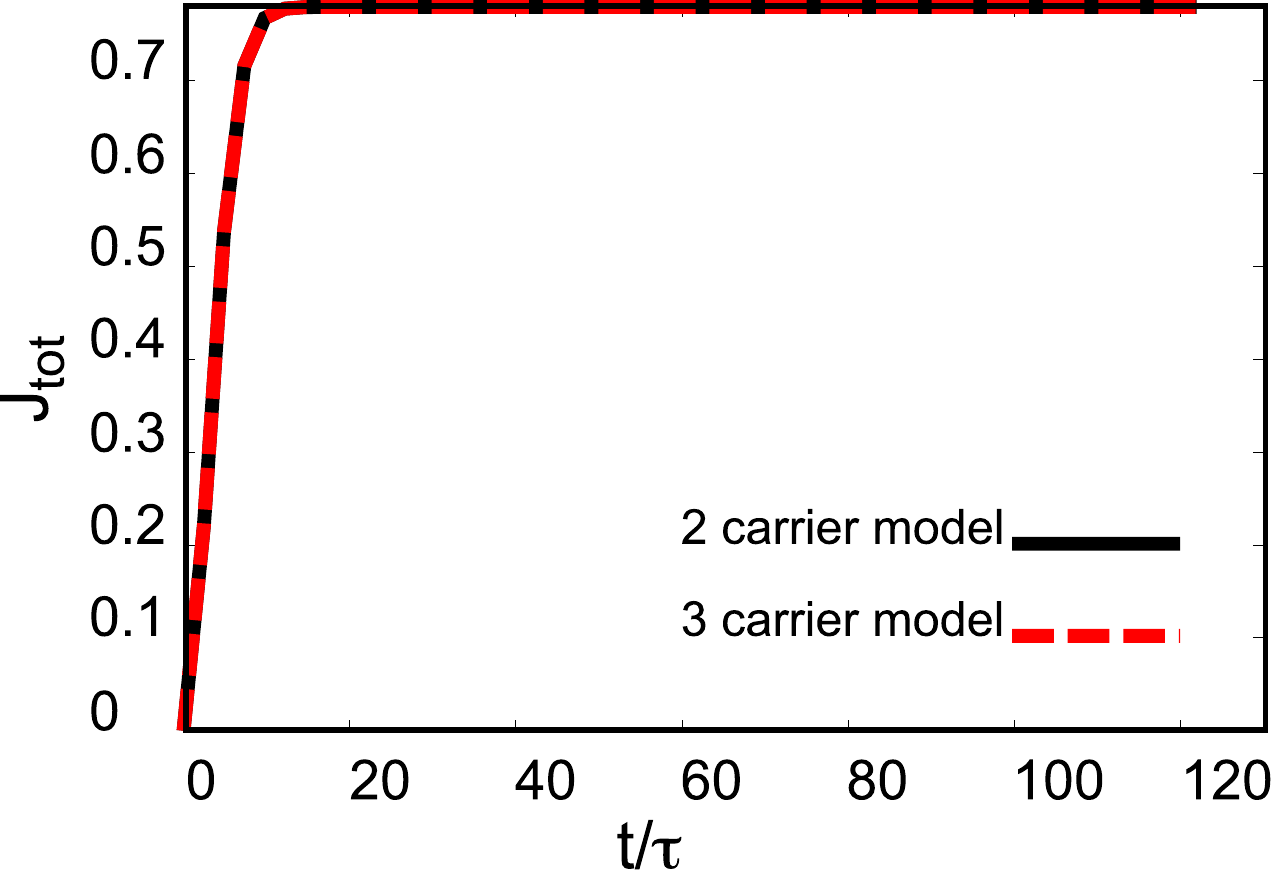}}
\caption{Photocurrent transient at low light intensity: effect 
of geminate pair recombination rate on rise time.}
\label{fig:G-kr}
\end{figure}

Figures~\ref{fig:G-u}-\ref{fig:G-kd}-\ref{fig:G-kr} refer to a device under 
low light intensity conditions and show the impact on the turn-on 
transient time of the value of the mobilities, of the geminate 
pair dissociation rate and of the recombination rate, respectively.

\begin{figure}[h]
\subfigure[Low mobilities: 
$\mu_{n}=\mu_{p} = 2\times 10^{-9} 
m^{2} V^{-1} s^{-1}$]{\includegraphics[width=.45\linewidth]{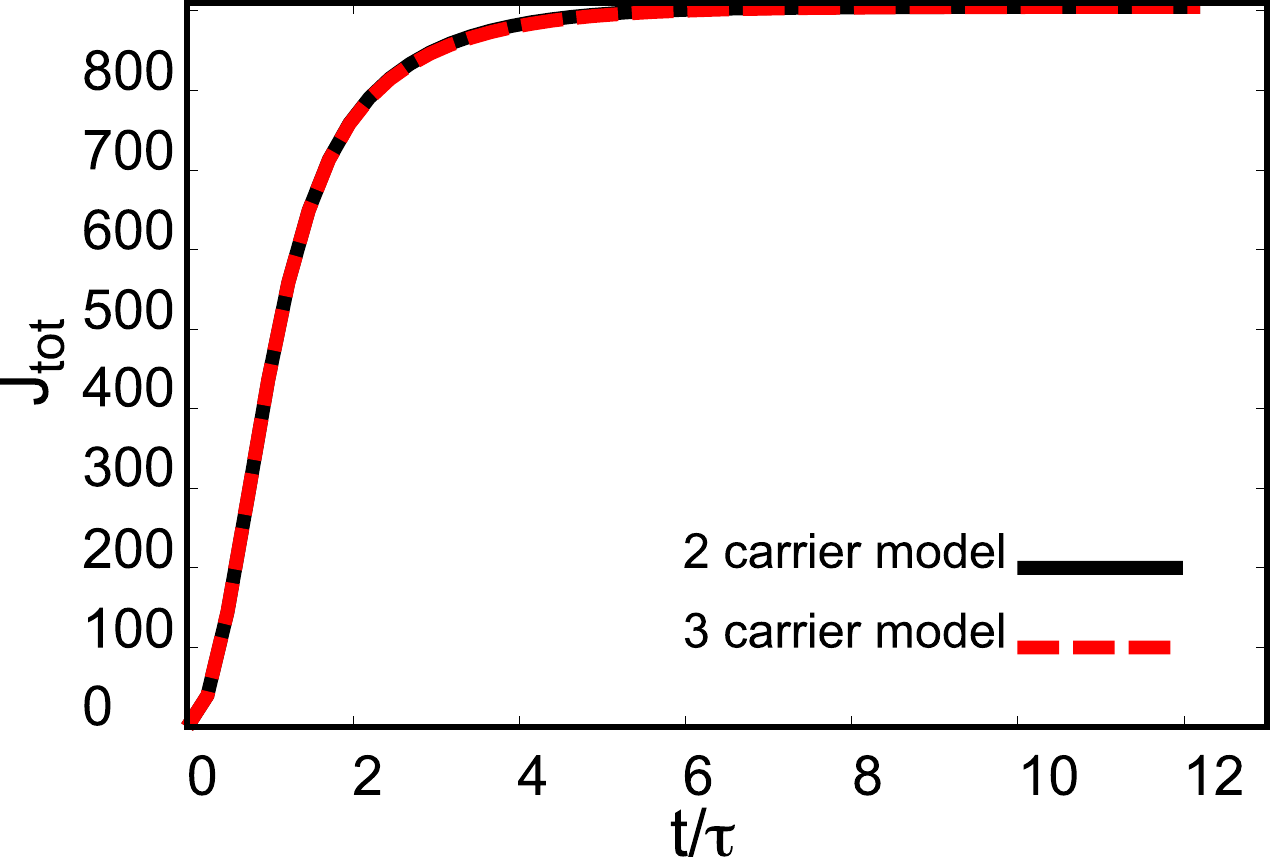}}
\subfigure[High mobilities: 
$\mu_{n}=\mu_{p} = 2\times 10^{-8} m^{2} V^{-1} s^{-1}$]
{\includegraphics[width=.45\linewidth]{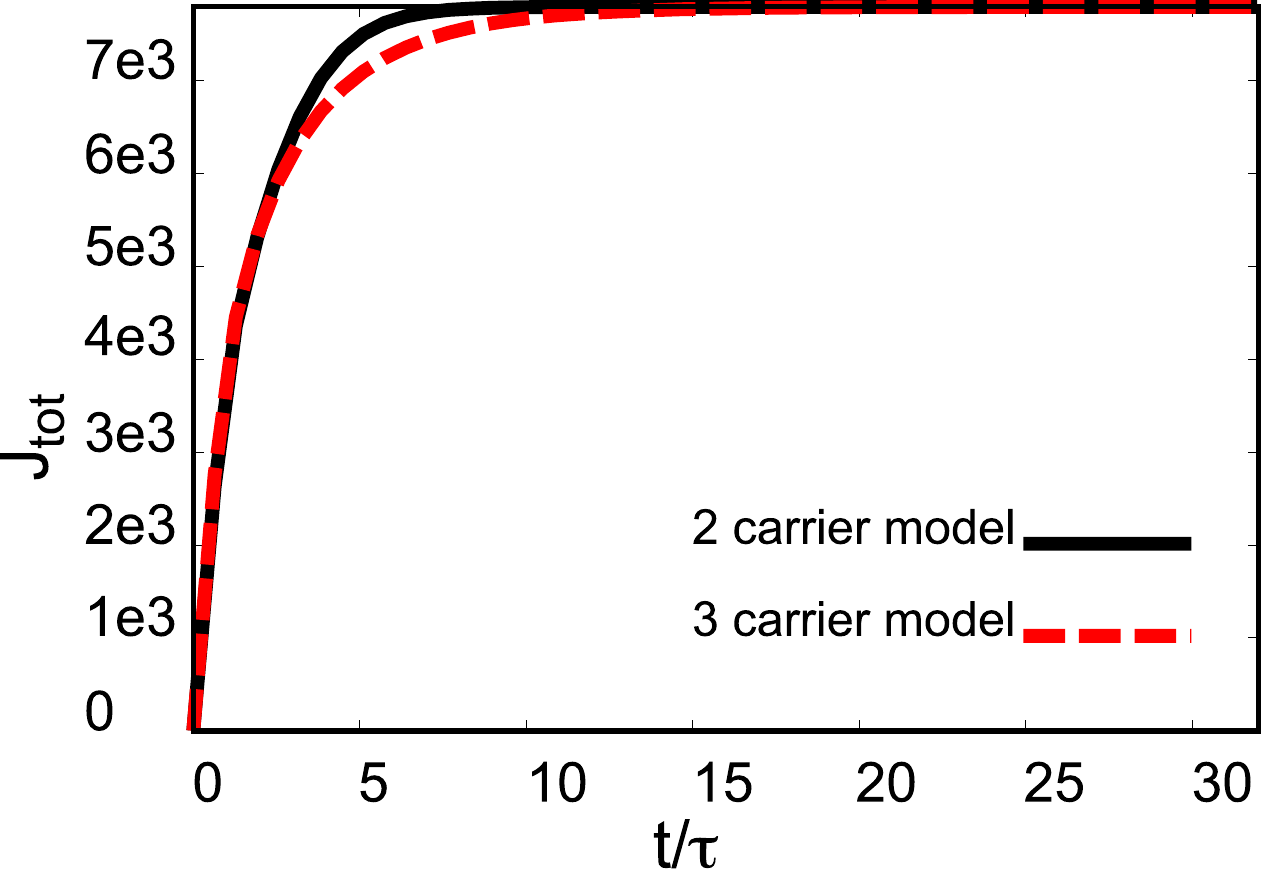}}
\caption{Photocurrent transient at high light intensity: 
effect of mobility on rise time.}
\label{fig:G+u}
\end{figure}

\begin{figure}[h]
\subfigure[$k_{diss}=4.4 \times 10^5 s^{-1}$]
{\includegraphics[width=.45\linewidth]{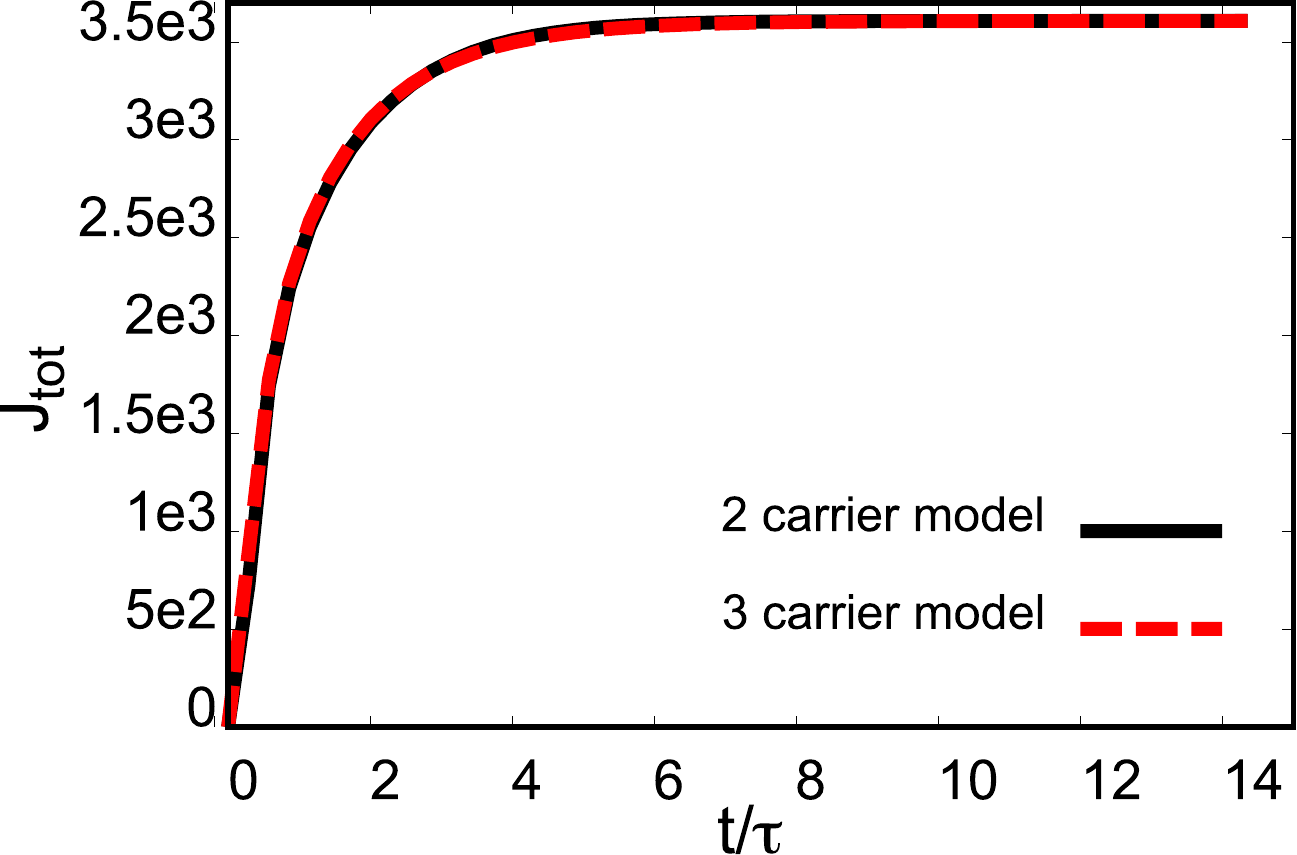}}
\subfigure[$k_{diss}=8 \times 10^6 s^{-1}$]
{\includegraphics[width=.45\linewidth]{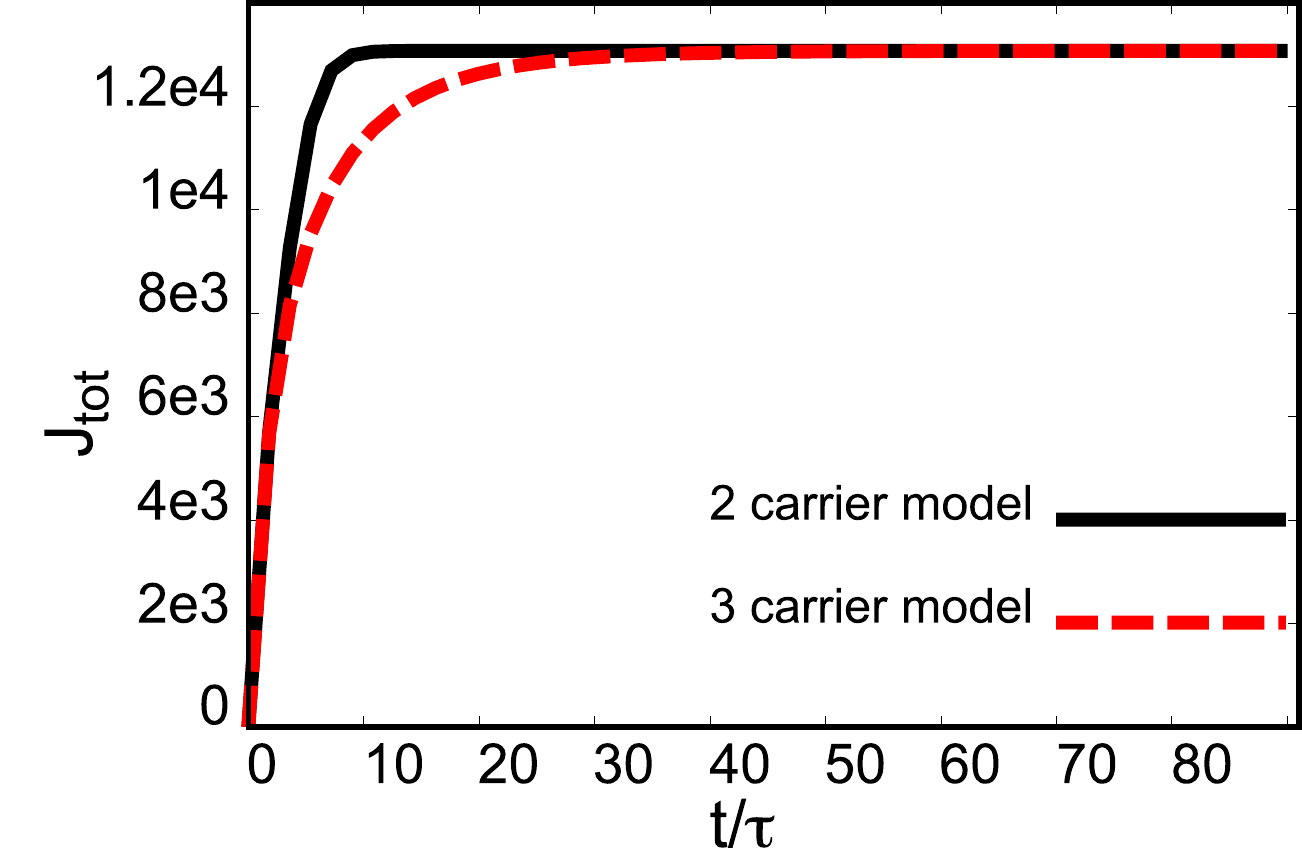}}
\caption{Photocurrent transient at high light intensity: 
effect of dissociation rate on rise time.}
\label{fig:G+kd}
\end{figure}

\begin{figure}[h]
\subfigure[$k_{rec}=10^5 s^{-1}$]
{\includegraphics[width=.45\linewidth]{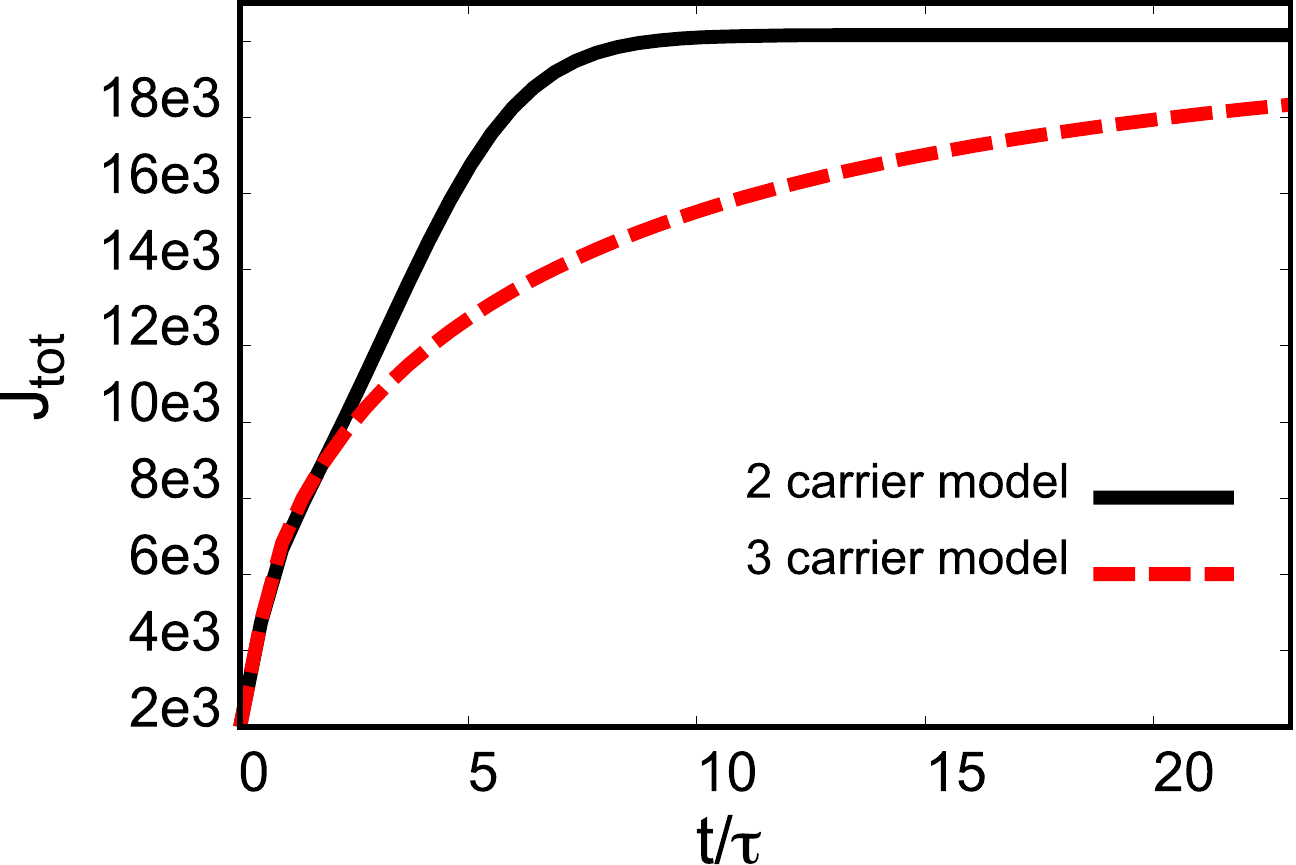}}
\subfigure[$k_{rec}=10^7 s^{-1}$]
{\includegraphics[width=.45\linewidth]{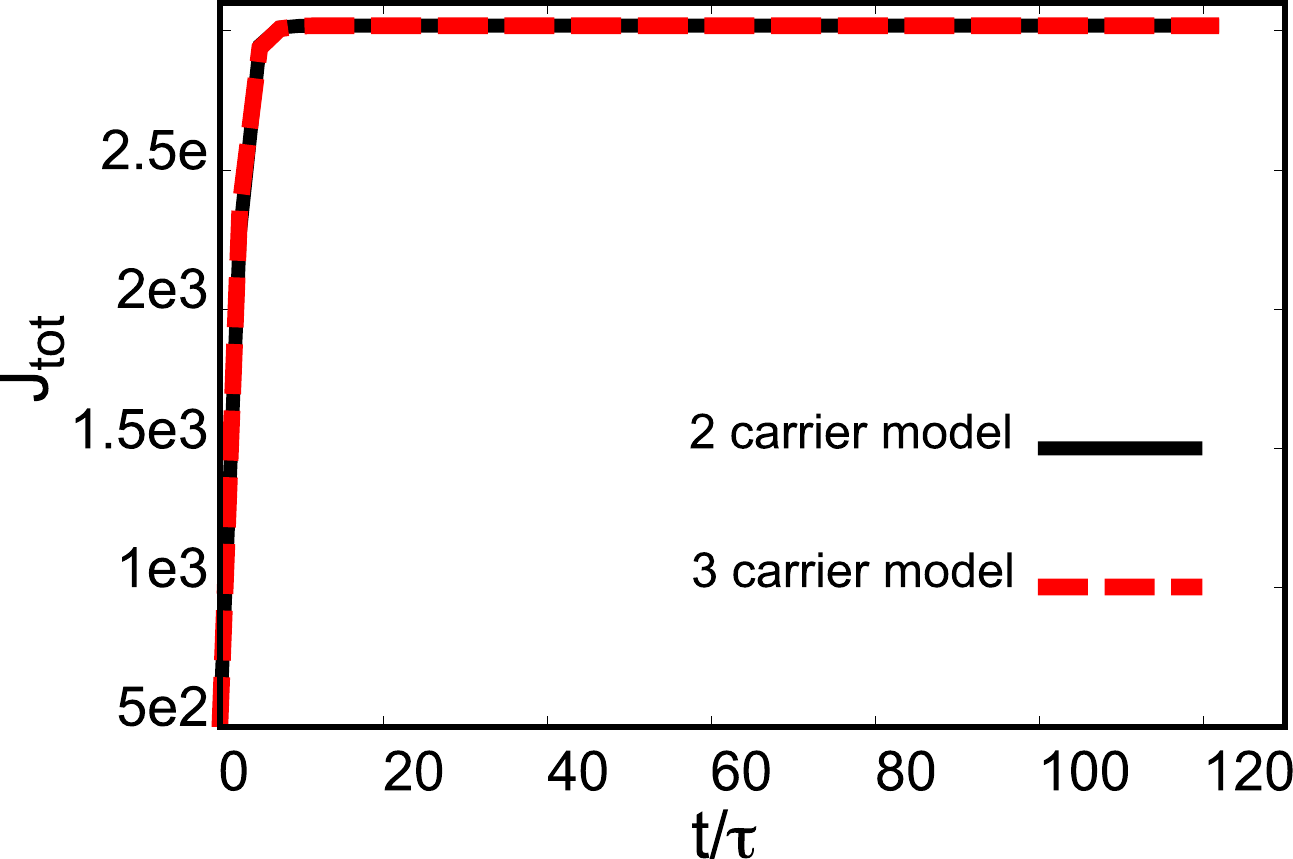}}
\caption{Photocurrent transient at high light intensity: 
effect of geminate pair recombination rate on rise time.}
\label{fig:G+kr}
\end{figure}

One may observe that, while at low intensity a change of one 
order of magnitude in the value of the mobility produces an almost 
equal change in the transient time, at high light intensity (Fig.~\ref{fig:G+u}) a similar 
change in the mobility has an almost negligible impact.
In this latter regime, variations in the dissociation rate 
$k_{diss}$ (Fig.~\ref{fig:G+kd}) and, more notably the recombination rate $k_{rec}$ (Fig.~\ref{fig:G+kr}), 
produce a more dramatic effect.

The analysis of the above results displays the complex relation 
between the transient 
behaviour of the device and the strongly nonlinear interplay among the several
occurring physical phenomena and shows the ability of the simplified 
model~\eqref{eq:model-transient}-\eqref{eq:eq_sys_const_rel} to capture such behaviour 
in most circumstances.
The only situation where the two models disagree is in the case of a device 
with high generation efficiency (i.e., a low value 
of $k_{rec}$) under high light intensity (cf.~Fig.~\ref{fig:G+kr}(a)).
Finally the steady-state current predicted by the reduced model is always in perfect agreement 
with that of the full model, as expected.

\section{Conclusions and Future Work}\label{sec:conclusions}

In this article, we have dealt with the mathematical
modeling and numerical simulation of photocurrent transients 
in nanostructured mono-layer OSCs. 
The model consists of a system of nonlinear diffusion-reaction PDEs with 
electrostatic convection, coupled to a kinetic ODE.
We have proposed a suitable reformulation of the model
which makes it similar to the drift-diffusion system for 
inorganic semiconductor devices.
This has allowed us to prove the existence of a solution for the problem 
in both stationary and transient conditions and to highlight 
the role of exciton dynamics in determining the device turn-on time. 
For the numerical treatment, we carried out a temporal 
semi-discretization using an implicit adaptive method, and 
the resulting sequence of differential subproblems was 
linearized using the Newton-Raphson method with inexact Jacobian. 
Exponentially fitted finite elements were used for spatial 
discretization, and a thorough validation of the computational 
model was carried out by extensively investigating 
the impact of the model parameters on photocurrent transient times.

Future work is warranted in the following three main areas: 
1) extensions to the model; 2) improvement of the analytical 
results; and 3) development of more specialized numerical algorithms. 
In detail:
\begin{description}
\item[1)] we intend to 
include exciton transport in order to be able to 
simulate multi-layer or 
nanostructured devices~\cite{Buxton2007,McGehee2007,advmater2006,nanostrctcells};
\item[2)]
we aim to extend Theorem~\ref{th:transient} to cover the full
problem~\eqref{eq:modelequations}--\eqref{eq:boundaryconditions}. 
A possible approach to achieve this result is to apply 
Theorem~\ref{th:transient} locally on a partition of $[0, T]$
into sub-intervals of size $\Delta t$, and verify the hypotheses 
of the Aubin lemma~\cite{Jerome1983} to extract a limiting solution as 
$\Delta t \rightarrow 0$;
\item[3)] starting from the above idea, we intend to devise 
a numerical algorithm for the local approximation of the full model system
over each sub-interval of size $\Delta t$ using the reduced 
model~\eqref{eq:model-transient}--\eqref{eq:eq_sys_const_rel}. 
The computer implementation of this approach is 
straightforward as it basically amounts to a successive application 
of the formulation discussed in Sect.~\ref{sec:numericalmethod} 
on each time slab. Furthermore, we intend to improve the robustness of the nonlinear solver 
with respect to the choice of scaling parameters 
(cf. Sect.~\ref{seq:spacedisc})
by adopting a staggered solution scheme
based on some variant of Gummel's Map~\cite{gummel64,ggmjcp}.
Such scheme could be either employed as an alternative to the current 
Newton solver or, even more effectively, combined with this latter
in a predictor-corrector fashion. The above modifications to the
solution algorithm are of great importance in dealing 
with the simulation of the multidimensional device structures 
mentioned at item 1).
\end{description}

\section{Acknowledgements} 
The very helpful comments by the anonymous referees contributed 
to improve the quality of the presentation, this contribution was much appreciated and is gratefully acknowledged.
The authors wish to thank Prof. Marco Sampietro and
Dr. Dario Natali, Dipartimento di Elettronica e Informazione,
Politecnico di Milano, Milano (Italy), 
for many stimulating discussions.
They also wish to thank Prof. Joseph W. Jerome for his very careful reading of 
the manuscript and for his useful comments and suggestions.
The first author
was partially supported by the European Research Council through 
the FP7 Ideas Starting Grant program {\it ``GeoPDEs -- Innovative compatible 
discretization techniques for Partial Differential Equations''}.
This support is gratefully acknowledged.

\bibliographystyle{elsart-num}
\bibliography{orgc}

\end{document}